\documentclass[12pt]{elsarticle}
\usepackage{amsfonts,enumerate,amsmath,amssymb,mathrsfs,amsbsy}
\usepackage[colorlinks=true]{hyperref}

\usepackage{caption}
\biboptions{numbers,sort&compress}

\usepackage{xcolor}

\def\qed{\strut\hfill $\Box$}
\newtheorem{thm}{Theorem}[section]

\newtheorem{lem}[thm]{Lemma}

\newtheorem{rem}[thm]{Remark}
\newtheorem{defn}[thm]{Definition}

\newcommand{\thmref}[1]{Theorem~{\rm \ref{#1}}}

\def\para#1{\vskip .4\baselineskip\noindent{\bf #1}}
\numberwithin{equation}{section}
\begin{document}
\begin{frontmatter}	
	\title{Averaging principles for non-autonomous two-time-scale stochastic reaction-diffusion equations with polynomial growth}
	
	\author[mymainaddress]{Ruifang Wang}
	\ead{wrfjy@yahoo.com}
	 
	\author[mymainaddress,myfivearyaddress]{Yong Xu\corref{mycorrespondingauthor}}
	\cortext[mycorrespondingauthor]{Corresponding author}
	\ead{hsux3@nwpu.edu.cn}
	
	\author[mysecondaryaddress]{Bin Pei}
	\ead{binpei@hotmail.com}
	
	
	
	\address[mymainaddress]{Department of Applied Mathematics, Northwestern Polytechnical University, Xi'an, 710072, China}
	\address[myfivearyaddress]{MIIT Key Laboratory of Dynamics and Control of Complex Systems, Northwestern Polytechnical University, Xi'an, 710072, China	}
	\address[mysecondaryaddress]{School of Mathematical Sciences, Fudan University, Shanghai,  200433, China}
	

	%

	\begin{abstract}
		In this paper, we develop the averaging principle for a class of two-time-scale stochastic reaction-diffusion equations driven by Wiener processes and Poisson random measures. We assume that all coefficients of the equation have polynomial growth, and the drift term of the equation is non-Lipschitz. Hence, the classical formulation of the averaging principle under the Lipschitz condition is no longer available. To prove the validity of the averaging principle, the existence and uniqueness of the mild solution
		are proved firstly. Then, the existence of time-dependent evolution family of measures associated with the fast equation is studied, by which the averaged coefficient is obtained. Finally, the validity of the averaging  principle is verified.	
	\end{abstract}		
\end{frontmatter}

\section{Introduction}\label{sec-1}
Multi-scale problems are widely encountered in composites, porous media, finance and other fields \cite{Harvey2011Multiple,Wu2016Approximate}. Morever, in practice, the parameters of systems often depend on time, 
non-autonomous systems are worthy of thorough analysis. 
For this reason, we are concerned with the following non-autonomous two-time-scale stochastic partial differential equations (SPDEs) on a bounded domain $ \mathcal{O} $  of  $ \mathbb{R}^d\left( d\ge 1 \right)  $: 
\small\begin{eqnarray}\label{orginal1}
\begin{aligned}
\begin{cases}
\frac{\partial u_{\epsilon}}{\partial t}\left( t,\xi \right) &=\mathcal{A}_1\left(t \right) u_{\epsilon}\left( t,\xi \right) +b_1\left(t, \xi ,u_{\epsilon}\left( t,\xi \right),v_{\epsilon}\left( t,\xi \right) \right) +f_1\left(t, \xi ,u_{\epsilon}\left( t,\xi \right) \right) \frac{\partial \omega ^{Q_1}}{\partial t}\left( t,\xi \right) \\
&\quad+\int_{\mathbb{Z}}{g_1\left(t, \xi ,u_{\epsilon}\left( t,\xi \right) ,z \right) \frac{\partial \tilde{N}_1}{\partial t}\left( t,\xi ,dz \right)},\\
\frac{\partial v_{\epsilon}}{\partial t}\left( t,\xi \right) 
&= \frac{1}{\epsilon}\left[ \left( \mathcal{A}_2\left( t \right) -\alpha \right) v_{\epsilon}\left( t,\xi \right) +b_2\left( t,\xi ,u_{\epsilon}\left( t,\xi \right) ,v_{\epsilon}\left( t,\xi \right) \right) \right] \\
&\quad+\frac{1}{\sqrt{\epsilon}}f_2\left( t, \xi, v_{\epsilon}\left( t,\xi \right) \right) \frac{\partial \omega ^{Q_2}}{\partial t}\left( t,\xi \right) +\int_{\mathbb{Z}}{g_2}\left( t, \xi, v_{\epsilon}\left( t,\xi \right) ,z \right) \frac{\partial \tilde{N}_{2}^{\epsilon}}{\partial t}\left( t,\xi ,dz \right) ,\\
u_{\epsilon}\left( 0,\xi \right) &=x\left( \xi \right), \quad 
v_{\epsilon}\left( 0,\xi \right) =y\left( \xi \right), \quad
\xi \in \mathcal{O},\cr
\mathcal{N}_1u_{\epsilon}\left( t,\xi \right) &= \mathcal{N}_2v_{\epsilon}\left( t,\xi \right) =0, \quad t\ge 0, \quad \xi \in \partial \mathcal{O},
\end{cases}
\end{aligned}
\end{eqnarray}
where $ \omega ^{Q_1}, \omega ^{Q_2}$ and $\tilde{N}_1, \tilde{N}_{2}^{\epsilon} $ are mutually independent Wiener processes and Poisson random measures,  $ 0<\epsilon \ll 1 $ is a  small parameter and $ \alpha $  is a sufficiently large fixed constant. In addition, $\mathcal{N}_i (i=1,2)$ are the boundary operators, which can be either the identity operator (Dirichlet boundary condition) or the first order operator (coefficients satisfying a uniform nontangentiality condition). The stochastic perturbations of the equations define on the same complete stochastic basis $\left( \varOmega ,\mathcal{F},\left\lbrace \mathcal{F}_t\right\rbrace _{t\geq0},\mathbb{P} \right) $, the specific introduction will be given in Section \ref{sec-2}.



The averaging principle is an effective method to analysis the slow-fast systems, which can simplify the system by constructing the averaged equation. In 1961, Bogolyubov and  Mitropolskii \cite{ Bogolyubov1961asymptotic}  studied the averaging principle, giving the first rigorous result for the deterministic case. Since then, the averaging principle became an active area of research. 
Khasminskii \cite{khas1968on} established the averaging principle for stochastic differential equations (SDEs) in 1968. Then, Givon \cite{givon2007strong}, Freidlin and Wentzell \cite{freidlin2012random},  Duan \cite{duan2014effective}, Xu and his co-workers \cite {xu2011averaging, xu2015approximation, xu2017stochastic} also studied the averaging principle of SDEs. In addition, many scholars also investigated the averaging principle of SPDEs in recent years, such as, Cerrai \cite{cerrai2009khasminskii,cerrai2011averaging}, Wang and Roberts \cite{wang2012average}, Pei and Xu \cite{pei2017two, pei2017stochastic, pei2017averaging}, Xu and Miao \cite{miao2017strong}. It should be pointed out that most of the current studies about the averaging principle are based on autonomous systems. In practical problems, the parameters of the system often depend on time. Therefore, non-autonomous system can depict some actual models better, which has made itself attract more and more attention of scholars. 

In 2017, effetive approximation for non-autonomous slow-fast system has been presented by Cerrai \cite{cerrai2017averaging},  and the system of this paper was driven by Gaussian noises.  
In our previous article \cite{Xu2018Averaging}, we study the non-autonomous slow-fast system driven by Gaussian noises and Poisson random measures.
An effective approximations for the slow equation of the original system in article \cite{Xu2018Averaging} was established by using the averaging principle, where the coefficients of the equation satisfy the Lipschitz condition and linear growth. But those conditions are too strict to study the validity of the averaging principle in many other relevant cases, such as, polynomial growth. One of the reaction-diffusion equations for the coefficients satisfy the polynomial growth is the Fitzhugh-Nagumo or Ginzburg-Landau type, those systems have appeared in the fields of biology and physics and attracted considerable attention. 
Therefore, we are devoted to developing the averaging principle for non-autonomous systems of reaction-diffusion equations with polynomial growth.

First, with the aid of the Sobolev embedding theorem, fixed point theorem and stopping technique, the existence and uniqueness of the mild solution is proved.
That is, for any $ p \geq 1 $ and $ T \geq s $,  we prove that system (\ref{orginal1}) admits a unique mild solution depending on the initial datum.

Next, as in our previous work \cite{Xu2018Averaging}, assuming the operator $  \mathcal{A}_2(t) $ is periodic and the functions $ b_2, f_2, g_2 $ are almost periodic.
Analyzing the fast equation with a frozen slow component and using Kunita’s
first inequality to deal with the Poisson terms,
we get that the evolution family of measures for the fast equation also exists, and it is almost periodic. Then, the averaged coefficient is defined through it, and the following averaged equation is obtained
\begin{eqnarray}\label{en15}
\begin{split}
\begin{cases}
\frac{\partial \bar{u}}{\partial t}\left( t,\xi \right) &=\mathcal{A}_1\left(t \right) \bar{u}\left( t,\xi \right) +\bar{b}_1\left(\xi, \bar{u}\left( t,\xi \right) \right) +f_1\left(t, \xi ,\bar{u}\left( t,\xi \right) \right) \frac{\partial \omega ^{Q_1}}{\partial t}\left( t,\xi \right) \\
&\quad+\int_{\mathbb{Z}}{g_1\left(t, \xi ,\bar{u}\left( t,\xi \right) ,z \right) \frac{\partial \tilde{N}_1}{\partial t}\left( t,\xi ,dz \right)},\\
\bar{u}\left( 0,\xi \right) &=x\left( \xi \right), \quad \xi\in \mathcal{O}, \quad
\mathcal{N}_1\bar{u}\left( t,\xi \right) =0, \quad t\ge 0,\quad \xi \in \partial \mathcal{O},
\end{cases} 
\end{split}
\end{eqnarray}
where $ \bar{b}_1\left(\xi, \bar{u}\left( t ,\xi\right) \right)  $ is the averaged coefficient,  which will be given in Section \ref{sec-5}.

Finally, the validity of the averaging principle is verified by using the classical Khasminskii method. That is, for any  $ T>0 $ and  $ \eta >0 $, we have 
\begin{eqnarray}\label{en711}
\underset{\epsilon \rightarrow 0}{\lim}\mathbb{P}\Big( \underset{t\in \left[ 0,T \right]}{\sup}\left\| u_{\epsilon}\left( t \right) -\bar{u}\left( t \right) \right\| _{\mathbb{D}(\bar{\mathcal{O}})}>\eta \Big) =0,
\end{eqnarray}
where  $ \bar{u} $ is the solution of the averaged equation (\ref{en15}).

We will give a specific definition of the notations in Section \ref{sec-2}. In this paper, $ c>0 $ with or without subscripts represents a general constant, the value of which may vary for different cases.
\section{Notations, assumptions and preliminaries}\label{sec-2}

Denote $ \mathbb{E} $ is the space $  \mathcal{C} (\bar{\mathcal{O}} ) $,  endowed with the following sup-norm 
$$
\left\| x \right\| _{ \mathbb{E}}=\underset{\xi \in \bar{\mathcal{O}}}{\text{sup}}\left| x\left( \xi \right) \right|,
$$
and the duality $ \left< \cdot ,\cdot \right> _{ \mathbb{E} } $. The norm of the product space $  \mathbb{E}\times\mathbb{E}$ denote as   
$$
\left\| x \right\| _{ \mathbb{E}\times  \mathbb{E}}=\left( \left\| x_1 \right\| _{ \mathbb{E}}^{2}+\left\| x_2 \right\| _{ \mathbb{E}}^{2} \right) ^{\frac{1}{2}}.
$$
and the corresponding duality of the product space $  \mathbb{E}\times\mathbb{E}$ is $ \left< \cdot ,\cdot \right> _{ \mathbb{E}\times \mathbb{E}} $. 

Let $ \mathbb{X} $ be any space, denote $\mathcal{L}\left( \mathbb{X} \right)$ is the space of the bounded linear operators in $\mathbb{X}$.  For any $ 0\leq s<T $ and  $ p\geq 1$, denote the norm of the space $ L^p\left( \Omega ;\mathbb{D}\left( \left[ s,T \right] ; \mathbb{X} \right)  \right) $ is
$$
\lVert u \rVert _{L_{s,T,p}\left( \mathbb{X} \right)}^{p}:=\mathbf{E}\underset{t\in \left[ s,T \right]}{\text{sup}}\lVert u\left( t \right) \rVert _{\mathbb{X}}^{p}.
$$
where $  \mathbb{D}\left( \left[ s,T \right] ; \mathbb{X} \right) $ denotes the space of all c\`adl\`ag path from $ \left[ s,T \right] $ into $ \mathbb{X} $.


For any $ p\in \left[1,\infty \right]  $ with $ p \neq 2 $, denote the norms of the space $ L^p\left(\mathcal{O} \right)  $ and $ L^p\left(\mathcal{O} \right)\times L^p\left(\mathcal{O} \right)  $ are both $ \left\| \cdot \right\| _p  $. When $ \delta > 0 $ and $ p<\infty $, we denote the norm of the space $ W^{\delta,p}\left(\mathcal{O}\right)  $ is $ \left\| \cdot\right\| _{\delta,p} $:
$$
\left\| x \right\| _{\delta ,p}=\left\| x \right\| _p+\Big( \int_D{\int_D{\frac{\left| x\left( \xi \right) -x\left( \eta \right) \right|}{\left| \xi -\eta \right|^{\delta p+d}}d\xi d\eta}} \Big) ^{\frac{1}{p}}.
$$

Now, we introduce some notations about subdifferential. 
The subdifferential of $ \left\| x \right\|_\mathbb{E} $ is defined as 
$$
\partial \left\| x \right\| _{ \mathbb{E}}:=\left\{ h\in  \mathbb{E}^*;\left\| h \right\| _{ \mathbb{E}^*}=1,\left< h,x \right> _{ \mathbb{E}}=\left\| x \right\| _{ \mathbb{E}} \right\},
$$
where $ \mathbb{E}^* $ is the dual space of $\mathbb{E}$.
Due to the characterization of the subdifferential \cite[Appendix D]{prato2014stochastic}, if $ u:\left[ 0,T\right] \rightarrow\mathbb{E} $  is any differentiable mapping, then
\begin{eqnarray}\label{en01}
\frac{d}{dt}^{-}\left\| u\left( t \right) \right\| _{ \mathbb{E}}\le \left< u{'}\left( t \right) ,\delta \right> _{ \mathbb{E}},
\end{eqnarray}
for any $ t\in \left[ 0,T\right]  $ and $ \delta \in \partial \left\| u\left( t \right) \right\| _{ \mathbb{E}} $.

Now, we assume the space dimension $d>1$, the processes ${\partial \omega ^{Q_1}}/{\partial t}\left( t,\xi \right)$ and ${\partial \omega ^{Q_2}}/{\partial t}\left( t,\xi \right)$ in the slow-fast system are the Gaussian noises, assumed it is white in time and colored in space. Here, $\omega ^{Q_i}\left( t,\xi \right) \left( i=1,2 \right)$ is the cylindrical Wiener processes, and it defined as $$
\omega ^{Q_i}\left( t,\xi \right) =\sum_{k=1}^{\infty}{Q_ie_k\left( \xi \right) \beta _k\left( t \right)}, \quad i=1,2,
$$
where $\left\{ e_k \right\} _{k\in \mathbb{N}}$ is a complete orthonormal basis in $ H $, $\left\{ \beta _k\left( t \right) \right\} _{k\in \mathbb{N}}$ is a sequence of mutually independent standard Brownian motion defined on the same complete stochastic basis $\left( \varOmega ,\mathcal{F},\left\lbrace \mathcal{F}_t\right\rbrace _{t\geq0},\mathbb{P} \right) $, and $Q_i$ is a bounded linear operator on $ H $.

Next, we give the definitions of Poisson random measures $\tilde{N}_1\left( dt,dz \right)$ and $\tilde{N}_{2}^{\epsilon}\left( dt,dz \right)$. Let $ \left( \mathbb{Z},\mathcal{B}\left( \mathbb{Z} \right) \right) $  be a given measurable space and  $ v\left( dz \right)  $ be a  $ \sigma  $-finite measure on it. $ D_{p_{t}^{i}},i=1,2 $ are two countable subsets of $ \mathbb{R}_+ $. Moreover, let
$ p_{t}^{1},t \in  D_{p_{t}^{1}}  $ be a stationary $ \mathcal{F}_t $-adapted Poisson point process on $ \mathbb{Z} $ with characteristic $ v $,
and $ p_{t}^{1},t \in  D_{p_{t}^{2}} $ be the other stationary $ \mathcal{F}_t $-adapted Poisson point process on $ \mathbb{Z} $ with characteristic $ {v}/{\epsilon} $. Denote by $ N_i\left( dt,dz \right),i=1,2 $ the Poisson counting measure associated with $ p_{t}^{i} $, i.e.,
$$
N_i\left( t,\varLambda \right) :=\sum_{s\in D_{p_{t}^{i}},s\le t}{I_{\varLambda}\left( p_{t}^{i} \right)},i=1,2.
$$
Let us denote the two independent compensated Poisson measures
$$\tilde{N}_1\left( dt,dz \right) :=N_1\left( dt,dz \right) -v_1\left( dz \right) dt$$ 
and 
$$\tilde{N}_{2}^{\epsilon}\left( dt,dz \right) :=N_2\left( dt,dz \right) -\frac{1}{\epsilon}v_2\left( dz \right) dt,$$ 
where  $v_1\left( dz \right) dt$ and $\frac{1}{\epsilon}v_2\left( dz \right) dt$ are the compensators.

In this paper, for any $ t\in \mathbb{R} $, the operators $\mathcal{A}_1\left( t\right) $ and the operators $\mathcal{A}_2\left( t\right) $ are the second order uniformly elliptic operators  with continuous coefficients on $\bar{\mathcal{O}}$. As in our previous work \cite{Xu2018Averaging}, we assume that the operator $ \mathcal{A}_i(t) $ has the following form
\begin{eqnarray}\label{en21}
\begin{split}
\mathcal{A}_i\left( t \right) =\gamma_i \left( t \right) \mathcal{A}_i+\mathcal{L}_i\left( t \right), \quad t\in \mathbb{R},\ i=1,2,
\end{split}
\end{eqnarray}
where $\mathcal{A}_i$ independent of $t$ is a second order uniformly elliptic operator, having continuous coefficients on $\bar{\mathcal{O}}$. And the operator $\mathcal{L}_i\left( t \right) $ is a first order differential operator has the following form
\begin{eqnarray}\label{en22}
\begin{split}
\mathcal{L}_i\left( t,\xi \right) u\left( \xi \right) =\left< l_i\left( t,\xi \right) ,\nabla u\left( \xi \right) \right> _{\mathbb{R}^d}, \quad t\in \mathbb{R}, \ \xi \in \bar{\mathcal{O}}.
\end{split}
\end{eqnarray}

Finally, for $ i=1,2 $, denote the realization of the operators $ \mathcal{A}_i $ and $ \mathcal{L}_i $  in $\mathbb{E}$  are $ A_i $ and $ L_i $, and the operator $ \mathcal{A}_i $ generates an analytic semigroup $ e^{tA_i} $.

Now, we give the following assumptions about the operators $ \mathcal{A}_i $ and $ Q_i $ as in \cite{Xu2018Averaging} and \cite{cerrai2011averaging}. 
\begin{enumerate}[({A}1)]
	\item
	\begin{enumerate}
		\item 	For $ i=1,2 $, the function $\gamma_i :\mathbb{R}\rightarrow \mathbb{R}$ is continuous, and there exist $\gamma _0, \gamma >0$ such that 
		\begin{eqnarray}\label{en23}
		\begin{split}
		\gamma _0\le \gamma_i \left( t \right) \le \gamma, \quad t\in \mathbb{R}.
		\end{split}
		\end{eqnarray}
		\item For $ i=1,2 $, the function $l_i:\mathbb{R}\times \bar{\mathcal{O}}\rightarrow \mathbb{R}^d$ is continuous and bounded.
	\end{enumerate}	
	\item For $i=1,2$, there exist a complete orthonormal system $\left\{ e_{i,k} \right\} _{k\in \mathbb{N}}$ of $ \mathbb{E}$, and two sequences of nonnegative real numbers $\left\{ \alpha _{i,k} \right\} _{k\in \mathbb{N}}$ and $\left\{ \lambda _{i,k} \right\} _{k\in \mathbb{N}}$ such that
	\begin{eqnarray}\label{en24}
	\begin{split}
	A_ie_{i,k}=-\alpha _{i,k}e_{i,k}, \quad Q_ie_{i,k}=\lambda _{i,k}e_{i,k}, \quad k\ge 1,
	\end{split}
	\end{eqnarray}
	and
	\begin{eqnarray}\label{en25}
	\begin{split}
	\kappa _i:=\sum_{k=1}^{\infty}{\lambda _{i,k}^{\rho _i}\left\| e_{i,k} \right\| _{\infty}^{2}}<\infty, \quad \zeta _i:=\sum_{k=1}^{\infty}{\alpha _{i,k}^{-\beta _i}\left\| e_{i,k} \right\| _{\infty}^{2}}<\infty,
	\end{split}
	\end{eqnarray}
	for some constants $\rho _i\in \left( 2,+\infty \right]$ and $\beta _i\in \left( 0,+\infty \right)$ such that
	\begin{eqnarray}\label{en26}
	[{\beta _i\left( \rho _i-2 \right)}]/{\rho _i}<1.
	\end{eqnarray}
\end{enumerate}	

About the coefficients of the system (\ref{orginal1}), we assume it satisfy the following conditions.
\begin{enumerate}[({A}3)]
	\item 
	\begin{enumerate}
		\item The mappings $b_1:\mathbb{R}\times\bar{\mathcal{O}}\times \mathbb{R}^2\rightarrow \mathbb{R} $ is continuous and there exists $ m_1\geq 1 $ such that
		\begin{eqnarray}\label{en1}
		\underset{\left(t, \xi\right)  \in  \mathbb{R}\times\bar{\mathcal{O}}  }{\text{sup}}\left| b_1\left( t,\xi ,x,y \right) \right|\le c\left( 1+\left| x \right|^{m_1}+\left| y \right| \right) ,\quad\left( x,y \right) \in \mathbb{R}^2.
		\end{eqnarray}
		\item There exists $c>0$ such that, for any $ x,h\in \mathbb{R}^2 $,
		\begin{eqnarray}\label{en04}
		\underset{\left(t, \xi\right)  \in   \mathbb{R}\times\bar{\mathcal{O}} }{\text{sup}}\left( b_1\left( t,\xi ,x+h \right) -b_1\left( t,\xi ,x \right) \right) h_1\le c\left| h_1 \right|\left( 1+\left| x \right|+\left| h \right| \right).
		\end{eqnarray}
		\item There exists $ \theta>0 $ such that
		\begin{eqnarray}\label{en03}
		\underset{\left(t, \xi\right)\in \mathbb{R}\times\bar{\mathcal{O}}}{\text{sup}}\left| b_1\left( t,\xi ,x \right) -b_1\left( t,\xi ,y \right) \right|\le c\left( 1+\left| x \right|^{\theta}+\left| y \right|^{\theta} \right) \left| x-y \right|,\quad x,y\in \mathbb{R}^2.
		\end{eqnarray}
	\end{enumerate}
\end{enumerate}	
\begin{enumerate}[({A}4)]
	\item 
	\begin{enumerate}
		\item The mappings $b_2:\mathbb{R}\times\bar{\mathcal{O}}\times \mathbb{R}^2\rightarrow \mathbb{R} $ is continuous and there exists $ m_2\geq 1 $ such that
		\begin{eqnarray}\label{en2}
		\underset{\left(t, \xi\right)  \in   \mathbb{R}\times\bar{\mathcal{O}} }{\text{sup}}\left| b_2\left( t,\xi ,x,y \right) \right|\le c\left( 1+\left| x \right|+\left| y \right|^{m_2} \right) ,\quad\left( x,y \right) \in \mathbb{R}^2.
		\end{eqnarray}
		\item There exists $c>0$ such that, for any $ x,h\in \mathbb{R}^2 $,
		\begin{eqnarray}\label{en05}
		\underset{\left(t, \xi\right)  \in   \mathbb{R}\times\bar{\mathcal{O}}}{\text{sup}}\left( b_2\left( t,\xi ,x+h \right) -b_2\left( t,\xi ,x \right) \right) h_2\le c\left| h_2 \right|\left( 1+\left| x \right|+\left| h \right| \right).
		\end{eqnarray}
		\item The mapping $ b_2\left(t,\xi,\cdot \right):\mathbb{R}^2\rightarrow\mathbb{R}  $ is locally Lipschitz-continuous, uniformly with respect to $ \left(t, \xi\right)  \in   \mathbb{R}\times\bar{\mathcal{O}} $. 
		\item	
		For all $x,y_1,y_2\in \mathbb{R}$, we have
		\begin{eqnarray}\label{en29}
		\begin{split}
		b_2\left( t,\xi ,x,y_1 \right) -b_2\left( t,\xi ,x,y_2 \right) =-\tau \left( t,\xi ,x,y_1,y_2 \right) \left( y_1-y_2 \right).
		\end{split}
		\end{eqnarray}
		for some measurable function $\tau :\mathbb{R}\times \bar{\mathcal{O}}\times \mathbb{R}^3\rightarrow \left[ 0,\infty \right)$.
	\end{enumerate}
\end{enumerate}
\begin{enumerate}[({A}5)]
	\item 
	The mappings $f_1:\mathbb{R}\times\bar{\mathcal{O}}\times \mathbb{R}\rightarrow \mathbb{R},g_1:\mathbb{R}\times\bar{\mathcal{O}}\times \mathbb{R}\times \mathbb{Z}\rightarrow \mathbb{R}, f_2:\mathbb{R}\times\bar{\mathcal{O}}\times \mathbb{R} \rightarrow \mathbb{R}, g_2:\mathbb{R}\times\bar{\mathcal{O}}\times \mathbb{R} \times \mathbb{Z}\rightarrow \mathbb{R}$ are continuous, and the mappings $f_1\left( t,\xi ,\cdot \right):\mathbb{R}\rightarrow \mathbb{R},g_1\left( t,\xi,\cdot,z \right):\mathbb{R}\rightarrow \mathbb{R},  f_2\left(t, \xi ,\cdot \right):\mathbb{R} \rightarrow \mathbb{R},g_2\left( t,\xi,\cdot,z \right):\mathbb{R} \rightarrow \mathbb{R}$ are Lipschitz-continuous, uniformly with respect to $(t,\xi,z) \in \mathbb{R}\times \bar{\mathcal{O}}\times \mathbb{Z}$. 
	Moreover, for all  $ p\ge 1 $, there exist positive constants $c_1,c_2 $,  such that for all $ x, y \in \mathbb{R} $, have
	$$
	\underset{\left( t,\xi \right) \in \mathbb{R}\times \bar{\mathcal{O}}}{\text{sup}}\int_{\mathbb{Z}}{\left|  g_i\left(t, \xi ,x,z \right) -g_i\left(t, \xi ,y,z \right)\right|  }^{p}\upsilon _i\left( dz \right) \leq c_i\left|  x-y\right| ^{p}, \quad i=1,2. 
	$$
\end{enumerate}
\begin{enumerate}[({A}6)]
	\item  For any $  x,y \in \mathbb{R}, $ it hold that
	\begin{eqnarray}
	\underset{\left( t,\xi \right) \in \mathbb{R}\times \bar{\mathcal{O}}}{\text{sup}}\Big(\left| f_i\left( t,\xi ,x \right) \right|^p+\int_{\mathbb{Z}}{\left| g_i\left( t,\xi ,x, z \right) \right|^p}v_i\left( dz \right) \Big)\le c\left( 1+\left| x \right|^{\frac{p}{m_i}} \right),\quad i=1,2,
	\end{eqnarray}
	where $ m_1 $ and $ m_2 $ are the constants introduced in (\ref{en1}) and (\ref{en2}).
\end{enumerate}
\begin{rem}\label{rem2.1}
	{\rm For any $\left( t,\xi \right) \in \mathbb{R}\times \bar{\mathcal{O}}$ and $x,y,h\in  \mathbb{E},z\in \mathbb{Z}$, we shall set
		\begin{gather}
		B_1\left(t, x,y \right) \left( \xi \right) :=b_1\left(t, \xi ,x\left( \xi \right) ,y\left( \xi \right) \right) , \quad B_2\left( t,x,y \right) \left( \xi \right) :=b_2\left( t,\xi ,x\left( \xi \right) ,y\left( \xi \right) \right) ,\cr
		\left[ F_1\left(t, x \right) h \right] \left( \xi \right) :=f_1\left(t, \xi ,x\left( \xi \right) \right) h\left( \xi \right), \quad \left[ F_2\left( t,x  \right) h \right] \left( \xi \right) :=f_2\left( t,\xi ,x\left( \xi \right)  \right) h\left( \xi \right) ,\cr
		\left[ G_1\left(t, x,z \right) h \right] \left( \xi \right) :=g_1\left(t, \xi ,x\left( \xi \right) ,z \right) h\left( \xi \right), \quad
		\left[ G_2\left( t,x,z \right) h \right] \left( \xi \right) :=g_2\left( t,\xi ,x\left( \xi \right),z \right) h\left( \xi \right).\nonumber
		\end{gather}
		Due to the assumption (A3)  and (A4), we know the mappings
		$ B_1 : \mathbb{R}\times \mathbb{E}\times \mathbb{E}\rightarrow \mathbb{E}$ and $ B_2: \mathbb{R}\times \mathbb{E}\times  \mathbb{E}\rightarrow \mathbb{E}$ are well defined and continuous. According to (\ref{en1}) and (\ref{en2}), for any $ x,y\in\mathbb{E}$ and $ t\in\mathbb{R} $, we have 
		\begin{eqnarray}\label{en02}
		\left\|  B_1\left( t,x,y \right) \right\|_ \mathbb{E}\le c\left( 1+\left\|  x \right\|_{ \mathbb{E}} ^{m_1}+\left\|  y \right\|_\mathbb{E} \right),\quad \left\|  B_2\left( t,x,y \right) \right\|_ \mathbb{E}\le c\left( 1+\left\|  x \right\|_{ \mathbb{E}} +\left\|  y \right\|_ \mathbb{E}^{m_2}  \right).
		\end{eqnarray}
		As a consequence of (\ref{en04}) and (\ref{en05}), it is immediate to check that, for any $ x,y,h,k\in\mathbb{E}$, any $ t\in\mathbb{R} $, and any $ \delta\in\partial \left\| h \right\| _{ \mathbb{E}} $, 
		\begin{eqnarray}\label{en06}
		\left< B_i\left( t,x+h,y+k \right) -B_i\left( t,x,y \right) ,\delta \right> _{ \mathbb{E}}\le c\left( 1+\left\| x \right\| _{ \mathbb{E}}+\left\| y \right\| _{ \mathbb{E}}+\left\| h \right\| _{ \mathbb{E}}+\left\| k \right\| _{ \mathbb{E}} \right) 
		\end{eqnarray}
		In view of (\ref{en03}), for any $x_1,y_1,x_2,y_2\in\mathbb{E}$, we have 
		\begin{small}
			\begin{align}\label{en0221}
			\!\left\| B_1\left( t,x_1,y_1 \right) -B_1\left( t,x_2,y_2 \right) \right\| _{\mathbb{E}}\le c\big( 1+\left\| \left( x_1,y_1 \right) \right\| _{\mathbb{E}\times \mathbb{E}}^{\theta}+\left\| \left( x_2,y_2 \right) \right\| _{\mathbb{E}\times \mathbb{E}}^{\theta} \big) \left( \left\| x_1-x_2 \right\| _{\mathbb{E}}+\left\| y_1-y_2 \right\| _{\mathbb{E}} \right).\!
			\end{align}
		\end{small}
		In addition, from the equation (\ref{en29}), for every $ \delta\in\partial \left\|k \right\| _{\mathbb{E}} $, we have 
		$$ \left\langle B_2\left( t,x,y+k \right)-B_2\left( t,x,y \right),\delta \right\rangle_\mathbb{E}\leq0  $$

		Due to  the assumption (A5) and  (A6), for any fixed $\left( t,z\right) \in\left( \mathbb{R},\mathbb{Z}\right) ,$ the mappings
		\begin{gather}
		F_i\left( t,\cdot \right) : \mathbb{E}\rightarrow \mathcal{L}\left( \mathbb{E}\right),\quad 
		G_i\left( t,\cdot ,z \right) : \mathbb{E}\rightarrow \mathcal{L}\left( \mathbb{E}\right),\quad i=1,2, \nonumber 
		\end{gather}
		are Lipschitz-continuous.}
\end{rem}

Now, for $ i=1,2 $, we define
$$
\gamma_i \left( t,s \right) :=\int_s^t{\gamma_i \left( r \right)}dr, \quad s<t,
$$
and let $ \gamma \left( t,s \right):=\left(\gamma_1 \left( t,s \right), \gamma_2 \left( t,s \right)  \right)   $. For any $\epsilon >0$ and $\beta \ge 0$, set
$$
U_{\beta ,\epsilon,i}\left( t,s \right) =e^{\frac{1}{\epsilon}\gamma_i \left( t,s \right) A_i-\frac{\beta}{\epsilon}\left( t-s \right)}, \quad s<t,
$$
in the case $\epsilon =1$, we write $U_{\beta,i}\left( t,s \right)$, and 
in the case $\epsilon =1$ and $\beta =0$, we write $U_i\left( t,s \right) $.

Next, for any $\epsilon >0,\beta \ge 0$ and for any $u\in  \mathbb{D} \left( \left[ s,t \right] ;\mathbb{E}\right) ,r\in \left[ s,t \right],$ we define
$$
\psi _{\beta ,\epsilon,i}\left( u;s \right) \left( r \right) =\frac{1}{\epsilon}\int_s^r{U_{\beta ,\epsilon,i}\left( r,\rho \right) L_i\left( \rho \right) u\left( \rho \right)}d\rho, \quad s<r<t,
$$
in the case $\epsilon =1$, we write $\psi _{\beta, i}\left( u;s \right) \left( r \right)$, and in the case $\epsilon =1$ and $\beta =0$, we write $\psi _{i}\left( u;s \right) \left( r \right)$.

We can easily get that $ \psi _{\beta ,\epsilon,i}\left( u;s \right) \left( t \right)  $ is the solution of
$$
du\left( t \right) =\frac{1}{\epsilon}\left( \mathcal{A}_i\left( t \right) -\beta \right) u\left( t \right) dt, \quad t>s, \ u\left( s \right) =0.
$$

\section{Existence, uniqueness of the solutions}\label{sec-3}
 More general, in this section, we mainly study the existence and uniqueness of solutions for the following problems 
\begin{eqnarray}\label{en032}
du\left( t \right)  = \left[   A(t)u\left( t \right) +B \left( t,u\left( t \right) \right) \right] dt+F \left( t,u\left( t \right) \right) dw^Q\left( t \right)  +\int_{\mathbb{Z}}{G\left( t,u\left( t \right) ,z \right) \tilde{N} \left( dt,dz \right)},
\end{eqnarray}
where
$$
u\left( t \right) :=\left( \begin{array}{c}
u_1\left( t \right)\\
u_2\left( t \right)\\
\end{array} \right) ,\quad A\left( t \right) :=\left( \begin{matrix}
A_1\left( t \right)&		0\\
0&		A_2\left( t \right) \\
\end{matrix} \right) ,\quad B \left( t,u\left( t \right) \right) =\left( \begin{array}{c}
B_1\left( t,u_1\left( t \right),u_2\left( t \right) \right)\\
B_2\left( t,u_1\left( t \right),u_2\left( t \right) \right)\\
\end{array} \right) ,
$$
and
$$\footnotesize
F \left( t,u\left( t \right)\right) :=\left( \begin{matrix}
F_1\left( t,u_1\left( t \right) \right)&		0\\
0&		 F_2\left( t,u_2\left( t \right) \right)\\
\end{matrix} \right) ,\quad G\left( t,u\left( t \right),z \right) :=\left( \begin{matrix}
G_1\left( t,u_1\left( t \right),z \right)&		0\\
0&		G_2\left( t,u_2\left( t \right),z \right)\\
\end{matrix} \right) ,
$$
and
$$
w^Q\left( t \right) :=\left( \begin{array}{c}
w^{Q_1}\left( t \right)\\
w^{Q_2}\left( t \right)\\
\end{array} \right) ,\quad \tilde{N} \left( dt,dz \right) :=\left( \begin{array}{c}
\tilde{N}_1\left( dt,dz \right)\\
\tilde{N}_{2}\left( dt,dz \right)\\
\end{array} \right) .\quad 
$$
According to the assumption (A5),  it is easy to know that for any fixed $\left( t,z\right) \in\left( \mathbb{R},\mathbb{Z}\right) ,$ the mappings
\begin{gather}
F\left( t,\cdot \right) : \mathbb{E}\times\mathbb{E}\rightarrow \mathcal{L}\left( \mathbb{E}\times\mathbb{E}\right),\quad 
G \left( t,\cdot ,z \right) : \mathbb{E}\times\mathbb{E}\rightarrow \mathcal{L}\left( \mathbb{E}\times\mathbb{E}\right), \nonumber 
\end{gather}
are Lipschitz-continuous.

\begin{defn}\label{defn3.1}
	For any fix $ (x_1,x_2)\in\mathbb{E}\times\mathbb{E}$, a process $ u(t) $ is a mild solution of the equation (\ref{en032}), if 
	\begin{eqnarray}\label{en034}
	u\left( t \right) &=&U \left( t,s \right)x+\psi \left( u;s \right)\left( t \right)+\int_s^t{U \left( t,r \right)B  \left(r, u \left( r \right)  \right)}dr\cr &&+\int_s^t{U \left( t,r \right) F  \left(r, u \left( r \right) \right)}dw^{Q }\left( r \right) \cr
	&&+\int_s^t{\int_{\mathbb{Z}}{U \left( t,r \right)G \left(r, u \left( r \right) ,z \right)}}\tilde{N}  \left( dr,dz \right), 
	\end{eqnarray}
\end{defn}
where
$$
U \left( t,s \right) =\left( \begin{matrix}
U_1\left( t,s \right)&		0\\
0&		U_{ 2}\left( t,s \right)\\
\end{matrix} \right), \  \psi \left( u ;s \right) \left( t \right) =\left( \begin{array}{c}
\psi _1\left( u_1;s \right) \left( t \right)\\
\psi _{2}\left( u_2;s \right) \left( t \right)\\
\end{array} \right), \   x=\left( \begin{array}{c}
x_1\\
x_2\\
\end{array} \right).
$$

Now, we denote
$$
\Phi(u)(t):=\int_s^t{U \left( t,r \right)B  \left(r, u \left( r \right)  \right)}dr,
$$
$$
\varGamma(u)(t):= \int_s^t{U \left( t,r \right)F  \left(r, u \left( r \right) \right)}dw^{Q }\left( r \right),
$$
and
$$
\varPsi(u)(t):=  \int_s^t{\int_{\mathbb{Z}}{U \left( t,r \right)G \left(r, u \left( r \right) ,z \right)}}\tilde{N}  \left( dr,dz \right).
$$

First, we prove that the mapping $ \varPsi(u)(t) $ is a contraction in $ L^p\left( \varOmega ;\mathbb{D}\left( \left[ s,T \right] ; \mathbb{E} \right) \right) $.
\begin{lem}\label{lem3.2}
	Under the assumptions (A1)-(A6), for any $ u,v\in L^p\left( \varOmega ;\mathbb{D}\left( \left[ s,T \right] ;\mathbb{E}\right) \right)  $ with $ p\geq 1 $,  the mapping $ \varPsi $ maps $ L^p\left( \varOmega ;\mathbb{D}\left( \left[ s,T \right] ;\mathbb{E}\right) \right)  $ into itself, and we have
	\begin{eqnarray}\label{en0}
	\left\|  \varPsi \left( u \right) -\varPsi \left( v \right) \right\| _{L_{s,T,p}(\mathbb{E})}\le c_{s,p}^{\varPsi}\left( T \right) \left\| u-v \right\|_{L_{s,T,p}(\mathbb{E})},
	\end{eqnarray}
	where $ c_{s,p}^{\varPsi} $ is a continuous increasing function with $ c_{s,p}^{\varPsi}\left(s \right)=0  $.
\end{lem} 
\para{Proof:} 
By using a factorization argument \cite[Theorem 8.3]{prato2014stochastic}, we have
$$
\varPsi \left( u \right) \left( t \right) -\varPsi \left( v \right) \left( t \right) =\frac{sin\pi \lambda}{\pi}\int_s^t{\left( t-r \right) ^{\lambda -1} U(t,r)\phi _{\lambda}\left( u,v \right) \left( r \right)}dr
$$
where
$$
\phi _{\lambda}\left( u,v \right) \left( r \right) :=\int_s^r{\int_{\mathbb{Z}}{\left( r-\sigma \right) ^{-\lambda}U(r,\sigma)\left[ G\left( \sigma ,u\left( \sigma \right) ,z \right) -G\left( \sigma ,v\left( \sigma \right) ,z \right) \right]}}\tilde{N} \left( d\sigma,dz \right),
$$
and $\lambda\in\left( 0,{1}/{2} \right) $. 

For any $ t,\epsilon>0 $  and $ p\geq1 $, the semigroup $ e^{tA} $ maps $ L^p(\mathcal{O};\mathbb{R}^2) $ into $ W^{\epsilon,p}(\mathcal{O};\mathbb{R}^2) $ and by using the semigroup law, we can obtain
\begin{eqnarray}\label{en021}
\left\| e^{tA}x \right\| _{\epsilon ,p}\le c \left( t\land 1 \right) ^{-\frac{\epsilon}{2}}\left\| x \right\| _p,\quad\ x\in L^p\left( \mathcal{O};\mathbb{R}^2 \right),
\end{eqnarray}
for some constant $ c $ independent of $ p $.
Then, according to (\ref{en021}), using the H\"{o}lder inequality, for any $ \epsilon < 2\lambda $, we have
\begin{eqnarray}\label{en033}
\lVert \varPsi \left( u \right) \left( t \right) -\varPsi \left( v \right) \left( t \right) \rVert _{\mathbb{E}}&\leq& 
\lVert \varPsi \left( u \right) \left( t \right) -\varPsi \left( v \right) \left( t \right) \rVert _{\epsilon ,p} \cr
&\le& c_{\lambda}\int_s^t{\left( \left( t-r \right) \land 1 \right) ^{\lambda -\frac{\epsilon}{2}-1}\lVert \phi _{\lambda}\left( u,v \right) \left( r \right) \rVert _p}dr\cr
&\le& c_{\lambda}\underset{r\in \left[ s,T \right]}{\text{sup}}\lVert \phi _{\lambda}\left( u,v \right) \left( r \right) \rVert _p\cdot\int_0^{t-s}{\left( r\land 1 \right) ^{\lambda -\frac{\epsilon}{2}-1}}dr,
\end{eqnarray}
so, if we show that $ \phi _{\lambda}\left( u,v \right) \left( r\right) \in L^p\left(  \mathcal{O};\mathbb{R}^2 \right)  $, we can get $ \varPsi(u)-\varPsi(v)\in  \mathbb{D} \left( \left[ s,T \right] ;W^{\epsilon ,p}\left( \mathcal{O};\mathbb{R}^2 \right) \right), $  $P-a.s. $
Using Kunita's first inequality \cite[Theorem 4.4.23]{Applebaum2009Processes} and the H\"{o}lder inequality, because $ G $ is Lipschitz-continuous, for any $ p\geq 1 $, we have
\begin{eqnarray}
\mathbf{E}\left| \phi _{\lambda}\left( u,v \right) \left( t,\xi \right) \right|^p 
&\le& c_p\mathbf{E}\Big( \int_s^t{\int_{\mathbb{Z}}{}}\left( t-\sigma \right) ^{-2\lambda}\left|  U(r,\sigma)\left[  G\left( \sigma ,u\left( \sigma \right) ,z \right) \left( \xi \right) \right. \right. \cr
&&\qquad\qquad\qquad\qquad\qquad\qquad\qquad-\left. \left.  G\left( \sigma ,v\left( \sigma \right) ,z \right) \left( \xi \right) \right] \right|^2v\left( dz \right) d\sigma \Big) ^{\frac{p}{2}}\cr
&&+c_p\mathbf{E}\int_s^t{\int_{\mathbb{Z}}{}}\left( t-\sigma \right) ^{-p\lambda}\left|  U(r,\sigma)\left[  G\left( \sigma ,u\left( \sigma \right) ,z \right) \left( \xi \right) \right.\right.  \cr
&&\qquad\qquad\qquad\qquad\qquad\qquad\qquad- \left. \left. G\left( \sigma ,v\left( \sigma \right) ,z \right) \left( \xi \right) \right] \right|^pv\left( dz \right) d\sigma\cr
&\le& c_p\mathbf{E}\Big( \int_s^t{\int_{\mathbb{Z}}{}}\left( t-\sigma \right) ^{-2\lambda}\left\|   U(r,\sigma)\left[   G\left( \sigma ,u\left( \sigma \right) ,z \right) \right. \right. \cr
&&\qquad\qquad\qquad\qquad\qquad\qquad\qquad-\left. \left.  G\left( \sigma ,v\left( \sigma \right) ,z \right) \right]  \right\| _\mathbb{E}^2v\left( dz \right) d\sigma \Big) ^{\frac{p}{2}}\cr
&&+c_p\mathbf{E}\int_s^t{\int_{\mathbb{Z}}{}}\left( t-\sigma \right) ^{-p\lambda}\left\|  U(r,\sigma)\left[    G\left( \sigma ,u\left( \sigma \right) ,z \right)  \right. \right. \cr
&&\qquad\qquad\qquad\qquad\qquad\qquad\qquad-\left. \left.  G\left( \sigma ,v\left( \sigma \right) ,z \right) \right]  \right\| _\mathbb{E}^pv\left( dz \right) d\sigma\cr
&\le& c_p\mathbf{E}\Big( \int_s^t{\left( t-\sigma \right) ^{-2\lambda}\left\|  u\left( \sigma   \right) -v\left( \sigma  \right)\right\|_\mathbb{E} ^2d\sigma} \Big) ^{\frac{p}{2}}\cr
&&+c_p\mathbf{E}\int_s^t{\left( t-\sigma \right) ^{-p\lambda}\left\|  u\left( \sigma   \right) -v\left( \sigma  \right)\right\|_\mathbb{E}^pd\sigma}\cr
&\le& c_{p,T} \underset{r\in \left[ s,T \right]}{\text{sup}}\mathbf{E}\lVert u\left( r \right) -v\left( r \right) \rVert _{\mathbb{E}}^{p}\cdot \Big[\Big( \int_0^{t-s}{ \sigma^{-\frac{2p\lambda}{p-2}}d\sigma}\Big) ^{\frac{p-2}{2}}+\int_0^{t-s}{ \sigma^{-p\lambda}d\sigma}\Big],\nonumber
\end{eqnarray}
so
\begin{eqnarray}\label{en036}
\mathbf{E}\lVert \phi _{\lambda}\left( u,v \right) \left( t \right) \rVert _{p} &=&\mathbf{E}\Big( \int_D{\left| \phi _{\lambda}\left( u,v \right) \left( t,\xi \right) \right|^p}d\xi \Big) ^{\frac{1}{p}}\cr
&\le& c_{p,T}|\mathcal{O}|^{\frac{1}{p}} \mathbf{E}\lVert u -v \rVert _{L_{s,T,p}(\mathbb{E})}\Big[\Big( \int_0^{t-s}{ \sigma^{-\frac{2p\lambda}{p-2}}d\sigma}\Big) ^{\frac{p-2}{2}}+\int_0^{t-s}{ \sigma^{-p\lambda}d\sigma}\Big]^{\frac{1}{p}}.
\end{eqnarray}
where $ |\mathcal{O}| $ is Lebesgue measure of the bounded domain $ \mathcal{O} $. Because $ u,v\in L^p\left( \varOmega ;\mathbb{D}\left( \left[ s,T \right] ;\mathbb{E}\right) \right)  $, so we know that $ \varPsi(u)-\varPsi(v)\in  \mathbb{D} \left( \left[ s,T \right] ;W^{\epsilon ,p}\left( \mathcal{O};\mathbb{R}^2 \right) \right), P-a.s. $ for any $ \lambda <\min(\frac{p-2}{2p},\frac{1}{p})    
$. 
In addition, we know
\begin{eqnarray}
\lVert \varPsi \left( u \right)  -\varPsi \left( v \right)  \rVert _{L_{s,T,p}\left( \mathbb{E} \right)} = \Big[ \underset{t\in \left[ s,T \right]}{\text{sup}}\lVert \varPsi \left( u \right) \left( t \right) -\varPsi \left( v \right) \left( t \right) \rVert _{\mathbb{E}}^p\Big] ^{\frac{1}{p}}
\end{eqnarray} 
according to the equation (\ref{en033}) and (\ref{en036}), we can get that $ \varPsi $ maps the space $ L^p\left( \varOmega ; \mathbb{D} \left( \left[ s,T \right] ;\mathbb{E}\right) \right)  $ into itself, and (\ref{en0}) holds with

\begin{eqnarray}
c_{s,p}^{\varPsi}\left( t \right) =\int_0^{t-s}{\left( r\land 1 \right) ^{\lambda -\frac{\epsilon}{2}-1}}dr\cdot\Big[\Big( \int_0^{t-s}{ \sigma^{-\frac{2p\lambda}{p-2}}d\sigma}\Big) ^{\frac{p-2}{2}}+\int_0^{t-s}{ \sigma^{-p\lambda}d\sigma}\Big]^{\frac{1}{p}}.\nonumber
\end{eqnarray}\qed
\begin{rem}\label{rem3.3}
	For any $ x:=\left( x_1,x_2\right) \in \mathbb{E}\times\mathbb{E}$, according to the assumption (A6), we know that there exists $ m:=\left(m_1,m_2 \right) \in \mathbb{R}\times\mathbb{R} $ and positive constants $ c $, such that 
	\begin{eqnarray}
	\underset{ \xi  \in \bar{\mathcal{O}}}{ \sup } \int_{\mathbb{Z}}{\left| G\left( t, x, z \right) \right|}v \left( dz \right) \le c\big( 1+\left| x \right|^{\frac{1}{m}} \big),\quad \left( t,z\right) \in\left( \mathbb{R},\mathbb{Z}\right). \nonumber
	\end{eqnarray}
	 For any $ p\geq 1 $, if $ u\in L^p\left( \varOmega ;\mathbb{D}\left( \left[ s,T \right] ;\mathbb{E} \right) \right)  $, by proceeding as Lemma \ref{lem3.2}, we can get that $ \varPsi \left( u \right)   \in  \mathbb{D}\left( \left[ s,T \right] ;W^{\epsilon ,p}\left(  {\mathcal{O}};\mathbb{R}^2 \right) \right) $, and it is easy to prove that there exists some continuous increasing function $ c_{s,p}^{\varPsi }\left( t\right) $ with $ c_{s,p}^{\varPsi }\left( s\right) =0 $, such that 
	\begin{eqnarray}
	\left\| \varPsi
	\left( u \right) \right\| _{L_{s,T,p}\left( \mathbb{E} \right)}^p\leq c_{s,p}^{\varPsi }\left( T \right) \big( 1+\left\| u \right\| _{L_{s,T,p}\left( \mathbb{E} \right)}^{\frac{p}{m}} \big).
	\end{eqnarray}
Moreover, as the space $ W^{\epsilon ,p}\left( \bar{\mathcal{O}};\mathbb{R}^2 \right)  $ continuously into $ C^{\theta} ( \bar{\mathcal{O}}  ) $ for any $ \theta<\epsilon-d/p $, so we have that 
$ \varPsi  \left( u \right) \in C^{\theta} ( \bar{\mathcal{O}}  ) $, and
\begin{eqnarray} \label{en5}
\mathbf{E}\underset{t\in \left[ s,T \right]}{\text{sup}}\left\| \varPsi  \left( u \right) \left( t \right) \right\| _{C^{\theta} ( \bar{\mathcal{O}}  )}^{p}\leq c_{s,p}^{\varPsi }\left( T \right) \Big( 1+\left\| u \right\| _{L_{s,T,p}\left(\mathbb{E}\right)}^{\frac{p}{m}} \Big).
\end{eqnarray} 		
\end{rem}

Now, for any $ \alpha>0 $ and $ u\in L^p\left( \varOmega ; \mathbb{D} \left( \left[ s,T \right] ;\mathbb{E} \right) \right)  $, we define 
$$
\varPsi _{\alpha}\left( u \right) \left( t \right) :=\int_s^t{\int_{\mathbb{Z}}{U_{\alpha}\left( t,r \right) G\left( r,u\left( r \right) ,z \right)}}\tilde{N}\left( dr,dz \right) .
$$
We also can prove that $ \varPsi_\alpha $ maps $ L^p\left( \varOmega ; \mathbb{D} \left( \left[ s,T \right] ;\mathbb{E} \right) \right) $ into itself for any $ p\geq 1 $ and
\begin{eqnarray}
\left\| \varPsi _{\alpha}\left( u \right) \right\|^p _{L_{s,T,p}\left( \mathbb{E} \right)}\leq c_{s,p}^{\varPsi ,\alpha}\left( T \right) \big( 1+\left\| u \right\| _{L_{s,T,p}\left( \mathbb{E} \right)}^{\frac{p}{m}} \big), 
\end{eqnarray}
for some continuous increasing function $ c_{s,p}^{\varPhi,\alpha}\left( t \right) $ and $ c_{s,p}^{\varPhi,\alpha}\left( s \right)=0 $.

Now, we prove the existence and uniqueness of the solution for system (\ref{en032}).
\begin{thm}\label{th3.4}
	Under the  assumptions (A1)-(A6), for any $ x\in\mathbb{E}$ and $ p\geq 1 $, there exists a unique mild solution $ u_s^x\in L^p\left( \Omega ;\mathbb{D}\left( \left[ s,T \right] ;\mathbb{E} \right)\right)  $ for equation (\ref{en032}). Moreover, there have 
	\begin{eqnarray}\label{en031}
	\left\|u  \right\| _{L_{s,T,p}\left(\mathbb{E}\right)}\le c_{s,p}\left( T \right) \left( 1+\left\| x \right\| _\mathbb{E} \right) ,
	\end{eqnarray}
	for some continuous increasing function $ c_{s,p} $.
\end{thm}
\para{Proof:} In order to prove the existence of the solution for  system (\ref{en032}), we construct the following equations. For any $ n\in \mathbb{N}, i=1,2 $ and $ (t,\xi)\in \left[0,\infty \right)\times\bar{\mathcal{O}}  $, we define 
$$
b_{i,n}\left( t,\xi ,\sigma \right) :=\left\{ \begin{array}{c}
b_i\left( t,\xi ,\sigma \right)\\
b_i\left( t,\xi ,n\sigma /\left| \sigma \right| \right)\\
\end{array} \right. \ \ \ \ \ \begin{array}{c}
if\ \left| \sigma \right|\leq n,\\
if\ \left| \sigma \right|>n.\\
\end{array}
$$
For any $ n\in \mathbb{N} $, we can easily know that $ b_{i,n}\left( t,\xi ,\cdot \right):\mathbb{R}\rightarrow\mathbb{R} $ is Lipschitz-continuous uniformly with respect to $ \xi\in \bar{\mathcal{O}} $ and $ t\in \left[ s,T\right]  $. For any $ x\in\mathbb{E}$, define the corresponding composition operator $ B_n $ associated with $ b_n=(b_{1,n}, b_{2,n}) $ is
$$
B_n(t,x)(\xi):=b_n(t,\xi,x(\xi)), \quad \xi\in \bar{\mathcal{O}}.
$$
It is easy to get that $ B_n(t,\cdot) $ is Lipschitz-continuous. Moreover, if $ m<n $, we have
\begin{eqnarray}\label{en039}
\left\| x \right\| _\mathbb{E}\leq m\Rightarrow B_m\left( t,x \right) =B_n\left( t,x \right) =B\left( t,x \right). 
\end{eqnarray}
Next, we give the following problem
\begin{eqnarray}\label{en0310}
du\left( t \right)  = \left[   A(t)u\left( t \right) +B_n\left( t,u\left( t \right) \right) \right] dt+F \left( t,u\left( t \right) \right) dw^Q\left( t \right)  +\int_{\mathbb{Z}}{G\left( t,u\left( t \right) ,z \right) \tilde{N} \left( dt,dz \right)}.
\end{eqnarray}
Because $ B_n(t,\cdot) $ is Lipschitz-continuous, so the mapping $\varPhi_n $
\begin{eqnarray}
\varPhi _n\left( u \right) \left( t \right) :=\int_s^t{U \left( t,r \right)B_n \left(r, u \left( r \right)  \right)}dr \nonumber
\end{eqnarray}
is Lipschitz-continuous in $ L^p\left( \varOmega ;\mathbb{D}\left( \left[s,T \right] ;\mathbb{E} \right)  \right)  $.  
By proceeding as \cite[Lemma 6.1.2]{CerraiSecond2001}, we can prove that for any $ t\in [s,T], \epsilon\in (0,1]$ and $ u,v \in L^p (\varOmega; \mathbb{D}([s,T];\mathbb{E})) $, it yield
\begin{eqnarray}\label{en0313}
\left\| \psi \left( u ;s \right) \left( t \right) \right\| _\mathbb{E}
&\leq&c\int_s^t{\left(  \left( t-r \right) \land 1 \right) ^{-\frac{1+\epsilon}{2}}\left\| u \left(r \right) \right\| _\mathbb{E}}dr\cr
&\leq&c\int_0^{t-s}{\left(   r   \land 1 \right) ^{-\frac{1+\epsilon}{2}}}dr\underset{r\in \left[ s,t \right]}{\text{sup}}\left\| u \left(r \right) \right\| _\mathbb{E},
\end{eqnarray}
so, for any $ p\geq 1,  $ we can get 
\begin{eqnarray}\label{en0318}
\left\| \psi \left( u   \right)  -\psi \left( v   \right)  \right\| _{L_{s,T,p}\left(\mathbb{E}\right) }
\leq\left\| \psi \left( u-v  \right) \right\| _{L_{s,T,p}\left(\mathbb{E}\right) }\leq c_{s,p}^{\psi}\left( T \right)\left\| u  -v \right\| _{L_{s,T,p}\left(\mathbb{E}\right)}.
\end{eqnarray} 
where $ c_{s,p}^{\psi } $ is a continuous increasing function  with $ c_{s,p}^{\psi }\left( s\right) =0 $. In addition, according to \cite[Theorem 4.2, Remark 4.3]{cerrai2003Stochastic}, we know that there exists a constant $ p_*\geq 1 $, such that for any $ p\geq p_* $, we have 
\begin{eqnarray}\label{en4}
\left\| \varGamma
\left( u \right)-\varGamma
\left( v \right) \right\| _{L_{s,T,p}\left( \mathbb{E} \right)}^p\leq c_{s,p}^{\varGamma }\left( T \right)  \left\| u-v \right\| _{L_{s,T,p}\left( \mathbb{E} \right)},
\end{eqnarray}
and
\begin{eqnarray}\label{en3}
\left\| \varGamma
\left( u \right) \right\| _{L_{s,T,p}\left( \mathbb{E} \right)}^p\leq c_{s,p}^{\varGamma }\left( T \right) \big( 1+\left\| u \right\| _{L_{s,T,p}\left( \mathbb{E} \right)}^{\frac{p}{m}} \big).
\end{eqnarray}
where $ c_{s,p}^{\varPsi } $ is a continuous increasing function  with $ c_{s,p}^{\varPsi }\left( s\right) =0 $.

Due to Lemma \ref{lem3.2}, we have get that the mapping $ \varPsi  $ is a contraction in $ L^p\left( \varOmega ;\mathbb{D}\left( \left[ s,T \right] ;\mathbb{E} \right) \right) $. Moreover, because $ \varPhi _n(u) $ is Lipschitz-continuous and according to the equation (\ref{en0318}) and (\ref{en4}), we can know that the mapping $  \varPhi _n , \psi $ and $ \varGamma $ are contraction in $  L^p\left( \varOmega ;\mathbb{D}\left(\left[  s,T \right] ;\mathbb{E} \right)\right)  $. So, we can get that the mild solution $ u_n $ of the equation (\ref{en0310}) is the unique fixed point of the following mapping 
$$
u\left( t \right) \mapsto U\left( t,s \right) x+\psi \left( u;s \right) \left( t \right) +\varPhi _n\left( u \right) \left( t \right) +\varGamma \left( u \right) \left( t \right) +\varPsi \left( u \right) \left( t \right). 
$$

Next, we prove that the sequence $ \left\lbrace u_n\right\rbrace  $ is bounded in $ L^p\left( \Omega ;\mathbb{D}\left( \left[  s,T \right] ;\mathbb{E} \right) \right). $
\begin{lem}
	For any $ n\in \mathbb{N} $ and $ t\in [s,T] $, there exists a continuous increaing function $ c_{s,p}(t) $ such that 
	\begin{eqnarray}\label{en0314}
	\left\| u_n \right\| _{L_{s,T,p}\left(\mathbb{E}\right)}\le c_{s,p}\left( T \right) \left( 1+\left\| x \right\| _\mathbb{E} \right) ,
	\end{eqnarray}
\end{lem}
\para{Proof:} Denote $ \Lambda(u_n) $ is the solution of 
\begin{eqnarray}\label{en0311}
dv\left( t \right)  = A(t)v\left( t \right)dt+F \left( t,u_n\left( t \right) \right) dw^Q\left( t \right)  +\int_{\mathbb{Z}}{G\left( t,u_n\left( t \right) ,z \right) \tilde{N} \left( dt,dz \right)},\quad  v(s)=0,
\end{eqnarray}
we can get that $ \Lambda(u_n)\in L^p\left( \varOmega ;\mathbb{D}\left( \left[ s,T \right] ;\mathbb{E} \right) \right)  $ is the unique fixed point of the mapping 
$$
v\left( t \right) \mapsto \psi \left( v;s \right) \left( t \right) +\varGamma \left( u_n \right) \left( t \right) +\varPsi \left( u_n \right) \left( t \right) .
$$
So, we have 
\begin{eqnarray}\label{en0312}
\left\| \varLambda \left( u_n \right) \left( t \right) \right\| _\mathbb{E}\le \left\| \psi \left( \varLambda \left( u_n \right) ;s \right) \left( t \right) \right\| _\mathbb{E}+\left\| \varGamma \left( u_n \right) \left( t \right) \right\| _\mathbb{E}+\left\| \varPsi \left( u_n \right) \left( t \right) \right\| _\mathbb{E}.
\end{eqnarray}
According to the equation (\ref{en0313}) and (\ref{en0312}), using the Gronwall inequality, we can get 
\begin{eqnarray}\label{en222}
\left\| \varLambda \left( u_n \right) \left( t \right) \right\|_\mathbb{E}&\le& c\int_s^t{\left(  \left( t-r \right) \land 1 \right) ^{-\frac{1+\epsilon}{2}}\left[ \left\| \varGamma \left( u_n \right) \left( r \right) \right\|_\mathbb{E}+\left\| \varPsi \left( u_n \right) \left( r \right) \right\|_\mathbb{E} \right] e^{c\int_r^{t}{\left(  \left( t -\sigma\right)  \land 1 \right) ^{-\frac{1+\epsilon}{2}}}d\sigma}}dr\cr
&&+\left\| \varGamma \left( u_n \right) \left( t \right) \right\|_\mathbb{E}+\left\| \varPsi \left( u_n \right) \left( t \right) \right\|_\mathbb{E}\cr
&\le& \Big( \underset{r\in \left[ s,T \right]}{\text{sup}}\left\| \varGamma \left( u_n \right) \left( r \right) \right\|_\mathbb{E}+\underset{r\in \left[ s,T \right]}{\text{sup}}\left\| \varPsi \left( u_n \right) \left( r \right) \right\|_\mathbb{E} \Big) \times \Big( e^{c\int_s^{t}{\left( \left( t-\sigma\right)  \land 1 \right) ^{-\frac{1+\epsilon}{2}}}d\sigma}-1 \Big) \cr
&&+\Big( \underset{r\in \left[ s,T \right]}{\text{sup}}\left\| \varGamma \left( u_n \right) \left( r \right) \right\|_\mathbb{E}+\underset{r\in \left[ s,T \right]}{\text{sup}}\left\| \varPsi \left( u_n \right) \left( r \right) \right\|_\mathbb{E} \Big) \cr
&\le& c_s\left( t \right) \Big( \underset{r\in \left[ s,T \right]}{\text{sup}}\left\| \varGamma \left( u_n \right) \left( r \right) \right\|_\mathbb{E}+\underset{r\in \left[ s,T \right]}{\text{sup}}\left\| \varPsi \left( u_n \right) \left( r \right) \right\|_\mathbb{E} \Big). 
\end{eqnarray}
If we set $ v_n\left( t \right) :=u_n\left( t \right) -\varLambda \left( u_n \right) \left( t \right)  $, we know that $ v_n $ is the solution of the problem
\begin{eqnarray}
\frac{dv_n}{dt}\left( t \right) =A\left( t \right) v_n\left( t \right) dt+B_n\left( t,v_n\left( t \right) +\varLambda \left( u_n \right) \left( t \right) \right) ,\quad v_n\left( s \right) =x.
\end{eqnarray} 
According to the assumptions (A3) and (A4), we know that there exists $ m:=\left(m_1,m_2 \right) \in \mathbb{R}\times\mathbb{R} $, such that for any $ \delta _{v_n}\in \partial \left\| v_n\left( t \right) \right\|_\mathbb{E} $, we yield
\begin{eqnarray}
\frac{d}{dt}^-\left\| v_n\left( t \right) \right\|_\mathbb{E}&\le& \left< A\left( t \right) v_n\left( t \right) ,\delta _{v_n} \right>_\mathbb{E}+\left< B_n\left( t,v_n\left( t \right) +\varLambda \left( u_n \right) \left( t \right) \right) ,\delta _{v_n} \right>_\mathbb{E}\cr
&\le& \left< A\left( t \right) v_n\left( t \right) ,\delta _{v_n} \right>_\mathbb{E}+\left< B_n\left( t,v_n\left( t \right) +\varLambda \left( u_n \right) \left( t \right) \right) \right. \cr
&&\left. \ -B_n\left( t,\varLambda \left( u_n \right) \left( t \right) \right) ,\delta _{v_n} \right>_\mathbb{E}+\left< B_n\left( t,\varLambda \left( u_n \right) \left( t \right) \right) ,\delta _{v_n} \right>_\mathbb{E}\cr
&\le& c\left\| v_n\left( t \right) \right\|_\mathbb{E}+c\left( 1+\left\| \varLambda \left( u_n \right) \left( t \right) \right\| _{\mathbb{E}}^{m} \right),\nonumber
\end{eqnarray}
so
\begin{eqnarray}\label{en111}
\left\| v_n\left( t \right) \right\|_\mathbb{E} &\le& e^{c\left( t-s \right)}\left\| x \right\|_\mathbb{E}+c\int_s^t{e^{c\left( t-r \right)}\left( 1+\left\| \varLambda \left( u_n \right) \left( r \right) \right\| _{\mathbb{E}}^{m} \right)}dr\cr
&\le&  c_s\left( t \right) \big(1+\left\| x \right\|_\mathbb{E}+ \underset{r\in \left[ s,T \right]}{\text{sup}}\left\| \varGamma \left( u_n \right) \left( r \right) \right\| _{\mathbb{E}}^{m}+\underset{r\in \left[ s,T \right]}{\text{sup}}\left\| \varPsi \left( u_n \right) \left( r \right) \right\| _{\mathbb{E}}^{m} \big).
\end{eqnarray}
Due to the definition of $ u_n(t) $ and the equation (\ref{en222}) and (\ref{en111}), we can get that
\begin{eqnarray}
\left\| u_n\left( t \right) \right\|_\mathbb{E}\le c_s\left( t \right) \Big( 1+\left\| x \right\|_\mathbb{E}+\underset{r\in \left[ s,T \right]}{\text{sup}}\left\| \varGamma \left( u_n \right) \left( r \right) \right\| _{\mathbb{E}}^{m}+\underset{r\in \left[ s,T \right]}{\text{sup}}\left\| \varPsi \left( u_n \right) \left( r \right) \right\| _{\mathbb{E}}^{m} \Big) .\nonumber
\end{eqnarray}
 So, due to (\ref{en3}) and Remark \ref{rem3.3}, we can get that there exists a constant $ p_*\geq 1 $, such that, for any $ p\geq p_* $, we have 
\begin{eqnarray}
\mathbf{E}\underset{r\in \left[ s,t \right]}{\text{sup}}\left\| u_n\left( r \right) \right\| _{\mathbb{E}}^{p}\le c_{s,p}\left( t \right) \Big( 1+\left\| x \right\|_\mathbb{E}^p+\left( c_{s,p}^{\varGamma}\left( t \right) +c_{s,p}^{\varPsi}\left( t \right) \right) \mathbf{E}\underset{r\in \left[ s,t \right]}{\text{sup}}\left\| u_n\left( r \right) \right\| _{\mathbb{E}}^{p} \Big) ,\nonumber
\end{eqnarray}
because $ c_{s,p}^{\varGamma}\left( s \right) = c_{s,p}^{\varPsi}\left( s \right) =0$ and $ c_{s,p}, c_{s,p}^{\varGamma}, c_{s,p}^{\varPsi} $ are continuous, there exists $ t_0 $, such that $ c_{s,p}\left( s+t_0 \right) \cdot[c_{s,p}^{\varGamma}\left( s+t_0 \right) +c_{s,p}^{\varPsi}\left( s+t_0 \right)]\leq 1/2 $. For any $ t \in [s,s+t_0] $ we have
\begin{eqnarray}
\mathbf{E}\underset{t\in \left[ s,s+t_0 \right]}{\text{sup}}\left\| u_n\left( t \right) \right\| _{\mathbb{E}}^{p}\le c_{s,p}\left( t \right) \left( 1+\left\| x \right\|_\mathbb{E}^p \right) .
\end{eqnarray}
By proceeding it in the intervals $ \left[ s+t_0,s+2t_0 \right] ,\left[ s+2t_0,s+3t_0 \right]  $ etc., we get that for any $ T>s $ and $ p\geq p_* $, (\ref{en0314}) holds. If $ p<p_* $, using  the H\"{o}lder inequality, we can get (\ref{en0314}) also holds. \qed

Finally, through the sequence $ \left\lbrace u_n^x\right\rbrace  $, we can prove that Theorem \ref{th3.4} holds. For any $ n\in \mathbb{N} $ and $ x\in\mathbb{E}$, we define
$$
\tau _{n} :=\text{inf}\left\{ t\geq s:\left\| u_n\left( t \right) \right\|_\mathbb{E}\ge n \right\},  
$$
and let
$$
\tau   :=\underset{n\in \mathbb{N}}{\text{sup}}\tau _{n} .
$$
We can prove that the sequence of stopping times $ \left\lbrace \tau_{n} \right\rbrace  $ is non-decreasing, and thanks to (\ref{en0314}), we can get that $ P(\tau   =+\infty)=1 $.

Therefore, for any $ t\geq s $ and $ w\in \left\lbrace \tau   =+\infty\right\rbrace  $, there exists $ m\in \mathbb{N} $ such that for any $ t\in[s,T] $, have $ t\leq \tau _{m} (w)$, and then we define
$$
u   (t)(w):=u_m  (t)(w).
$$
Set $ \eta:=\tau _{n} \land\tau _{m}  $, due to (\ref{en039}), we can get
\begin{eqnarray}\label{en0316}
\left\| u_m  \left( t\land \eta \right) -u_n  \left( t\land \eta \right) \right\|_\mathbb{E}
&=&\left\| \psi \left( u_m  -u_n  ;s \right) \left( t\land \eta \right) \right\|_\mathbb{E}\cr
&&+\Big\| \int_s^{t\land \eta}{U\left( t\land \eta ,r \right) \left[ B_m\left( r,u_m  \left( r \right) \right) -B_n\left( r,u_n  \left( r \right) \right) \right]}dr \Big\|_\mathbb{E}\cr
&&+\left\| \varGamma \left( u_m   \right) \left( t\land \eta \right) -\varGamma \left( u_n   \right) \left( t\land \eta \right) \right\|_\mathbb{E}\cr
&&+\left\| \varPsi \left( u_m   \right) \left( t\land \eta \right) -\varPsi \left( u_n   \right) \left( t\land \eta \right) \right\|_\mathbb{E}\cr
&=&\left\| \psi \left( u_m  -u_n  ;s \right) \left( t\land \eta \right) \right\|_\mathbb{E}+\Big\| \int_s^t{I_{\left\{ r\le \eta \right\}}U\left( t\land \eta ,r \right)} \cr
&&\quad\times\left[ B_{m\vee n}\left( r\land \eta ,u_m  \left( r\land \eta \right) \right) -B_{m\vee n}\left( r\land \eta ,u_n  \left( r\land \eta \right) \right) \right]dr \Big\|_\mathbb{E}\cr
&&+\left\| \varGamma \left( u_m   \right) \left( t\land \eta \right) -\varGamma \left( u_n   \right) \left( t\land \eta \right) \right\|_\mathbb{E}\cr
&&+\left\| \varPsi \left( u_m   \right) \left( t\land \eta \right) -\varPsi \left( u_n   \right) \left( t\land \eta \right) \right\|_\mathbb{E}\cr
&\le& \underset{r\in \left[ s,t \right]}{\text{sup}}\left\| \psi \left( u_m  -u_n  ;s \right) \left( r\land \eta \right) \right\|_\mathbb{E}\cr
&&+c\int_s^t{\underset{\sigma \in \left[ s,r \right]}{\text{sup}}\left\| u_m  \left( \sigma \land \eta \right) -u_n  \left( \sigma \land \eta \right) \right\|_\mathbb{E}}dr\cr
&&+\underset{r\in \left[ s,t \right]}{\text{sup}}\left\| \varGamma \left( u_m   \right) \left( r\land \eta \right) -\varGamma \left( u_n   \right) \left( r\land \eta \right) \right\|_\mathbb{E}\cr
&&+\underset{r\in \left[ s,t \right]}{\text{sup}}\left\| \varPsi \left( u_m   \right) \left( r\land \eta \right) -\varPsi \left( u_n   \right) \left( r\land \eta \right) \right\|_\mathbb{E} . 
\end{eqnarray}
By proceeding as the proof of Lemma \ref{lem3.2}, using the factorization arguement for $ \varPsi \left( u_m   \right) \left( r\land \eta \right) -\varPsi \left( u_n   \right) \left( r\land \eta \right)   $,  we can obtain 
\begin{eqnarray}\label{en0317}
\mathbf{E}\underset{r\in \left[ s,t \right]}{\text{sup}}\left\| \varPsi \left( u_m   \right) \left( r\land \eta \right) -\varPsi \left( u_n   \right) \left( r\land \eta \right)  \right\|_\mathbb{E}\leq c_{s,1}^{\varPsi}\left( t \right) \mathbf{E}\underset{r\in \left[ s,t \right]}{\text{sup}}\left\|  u_m\left( r\land \eta \right)  -u_n\left( r\land \eta \right) \right\|_\mathbb{E}.
\end{eqnarray}
Then, substitue (\ref{en0318}), (\ref{en4}) and  (\ref{en0317}) into (\ref{en0316}), we have
\begin{eqnarray}
\underset{r\in \left[ s,t \right]}{\text{sup}}\left\|  u_m\left( r\land \eta \right)  -u_n\left( r\land \eta \right) \right\|_\mathbb{E}&\leq& \big( c_{s,1}^{\psi}\left( t \right) +c_{s,1}^{\varGamma}\left( t \right) +c_{s,1}^{\varPsi}\left( t \right) \big) \underset{r\in \left[ s,t \right]}{\text{sup}}\left\| u_m\left( r\land \eta \right)  -u_n\left( r\land \eta \right) \right\|_\mathbb{E}\cr
&&+c\int_s^t{\underset{\sigma \in \left[ s,r \right]}{\text{sup}}\left\| u_m\left( r\land \eta \right)  -u_n\left( r\land \eta \right) \right\|_\mathbb{E}}dr.\nonumber
\end{eqnarray}
Fix $ t_0>0 $, such that $ c_{s,1}^{\psi}\left( t_0 \right) +c_{s,1}^{\varGamma}\left( t_0 \right) +c_{s,1}^{\varPsi}\left( t_0 \right) \leq 1/2 $, we can get 
\begin{eqnarray}
\mathbf{E}\underset{r\in \left[ s,s+t_0 \right]}{\text{sup}}\left\|  u_m\left( r\land \eta \right)  -u_n\left( r\land \eta \right) \right\|_\mathbb{E}\leq c\int_s^{s+t_0}{\mathbf{E}\underset{\sigma \in \left[ s,r \right]}{\text{sup}}\left\|  u_m\left( r\land \eta \right)  -u_n\left( r\land \eta \right) \right\|_\mathbb{E}}dr.\nonumber
\end{eqnarray}
According to the Gronwall lemma, we have $ \mathbf{E}\underset{r\in \left[ s,s+t_0 \right]}{\text{sup}}\left\|  u_m\left( r\land \eta \right)  -u_n\left( r\land \eta \right) \right\|_\mathbb{E}=0 $, that is, for any $ t\in [s,s+t_0] $, we have $ u_m  (t\land\eta)=u_n  (t\land\eta). $ Repeat it in the interval $ [s+t_0,s+2t_0], [s+2t_0,s+3t_0],$ etc., we obtain
\begin{eqnarray}\label{en320}
u_m  \left( t \right) =u_n  \left( t \right) ,\quad s\leq t\leq \tau _{m} \land \tau _{n}, 
\end{eqnarray}
for any $ n\in \mathbb{N} $.
Because when $ w\in \left\lbrace \tau  =+\infty \right\rbrace  $ and $ t\leq \tau_m   $, we have denote $ u     (t)=u_m  (t) $, thanks to (\ref{en039}), this yields
\begin{eqnarray}
u     \left( t \right) &=&U \left( t,s\right)x+\psi \left( u     ;s \right)\left( t \right)+\int  _s ^t{U \left( t,r \right)B \left(r, u      \left( r \right)  \right)}dr +\int  _s ^t{U \left( t,r \right)F  \left(r, u      \left( r \right) \right)}dw^{Q }\left( r \right) \cr
&&+\int   _s^t{\int_{\mathbb{Z}}{U \left( t,r \right)G \left(r, u      \left( r \right) ,z \right)}}\tilde{N}  \left( dr,dz \right), \nonumber
\end{eqnarray}
$ P-a.s. $, that is, $ u     (t) $ is the mild solution of the system (\ref{en032}).

Now, we prove the solution of system (\ref{en032}) is unique. Denote another solution of system (\ref{en032}) is $ v    $, by proceeding as the equation (\ref{en320}), we can get that for any $ n\in \mathbb{N} $
\begin{eqnarray}
u  \left( t     \right) =v  \left( t  \right) ,\quad s\leq   t\leq \tau _{n}.\nonumber
\end{eqnarray}
For any $ T\geq s $, we know $ \left\lbrace \tau_n  \leq T \right\rbrace \downarrow \left\lbrace \tau  \leq T \right\rbrace$, we get that $ u     =v      $.

Finally, for any $ p\geq 1 $ and $ T>s $, we have 
\begin{eqnarray}
\underset{t\in \left[ s,T \right]}{\text{sup}}\lVert u  \left( t \right) \rVert _{\mathbb{E}}^{p}=\underset{n\rightarrow +\infty}{\lim}\underset{t\in \left[ s,T \right]}{\text{sup}}\lVert u  \left( t \right) \rVert _{\mathbb{E}}^{p}I_{\left\{ T\leq \tau _{n}  \right\}}=\underset{n\rightarrow +\infty}{\lim}\underset{t\in \left[ s,T \right]}{\text{sup}}\lVert u_{n} \left( t \right) \rVert _{\mathbb{E}}^{p}I_{\left\{ T\leq \tau _{n}  \right\}},\nonumber
\end{eqnarray}
according  to the estimate (\ref{en0314}) and the Fatou lemma, we can get (\ref{en031}). 

\section{The slow-fast system}\label{sec-4}
According to the introduced in Section \ref{sec-2}, system (\ref{orginal1}) can be rewritten as:
\begin{eqnarray}\label{orginal2}
\begin{split}
\begin{cases}
du_{\epsilon}\left( t \right) &=\left[ A_1\left( t\right) u_{\epsilon}\left( t \right) +B_1\left(t, u_{\epsilon}\left( t \right) ,v_{\epsilon}\left( t \right) \right) \right] dt+F_1\left(t, u_{\epsilon}\left( t \right) \right) d\omega ^{Q_1}\left( t \right) \\
&\quad +\int_{\mathbb{Z}}{G_1\left(t, u_{\epsilon}\left( t \right) ,z \right)}\tilde{N}_1\left( dt,dz \right), \\
dv_{\epsilon}\left( t \right) &=\frac{1}{\epsilon}\left[ \left( A_2\left( t \right) -\alpha \right) v_{\epsilon}\left( t \right) +B_2\left( t,u_{\epsilon}\left( t \right) ,v_{\epsilon}\left( t \right) \right) \right] dt \\
&\quad +\frac{1}{\sqrt{\epsilon}}F_2\left( t,v_{\epsilon}\left( t \right) \right) d\omega ^{Q_2}\left( t \right)  +\int_{\mathbb{Z}}{G_2}\left( t,v_{\epsilon}\left( t \right)  ,z\right) \tilde{N}_{2}^{\epsilon}\left( dt,dz \right), \\
u_{\epsilon}\left( s \right) &=x, \quad v_{\epsilon}\left( s \right) =y.
\end{cases}
\end{split}
\end{eqnarray}

Since the coefficients under the  assumptions (A1)-(A6) are uniform with respect to $ t\in \mathbb{R} $, according  the  prove in Section \ref{sec-3}, we can get that there exist two unique adapted $ u_{\epsilon} $ and $ v_{\epsilon} $ in $ L^p\left( \Omega ;\mathbb{D}\left( \left[ s,T \right] ;\mathbb{E} \right)  \right) $, such that
\begin{eqnarray}\label{or32}
\begin{split}
\begin{cases}
u_{\epsilon}\left( t \right) &=U_{1 }\left( t,s\right)x+\psi _{1}\left( u_{\epsilon};s \right)\left( t \right)+\int_s^t{U_{1}\left( t,r \right)B_1\left(r, u_{\epsilon}\left( r \right) ,v_{\epsilon}\left( r \right) \right)}dr\cr
&\quad+\int_s^t{U_{1}\left( t,r \right)F_1\left(r, u_{\epsilon}\left( r \right) \right)}dw^{Q_1}\left( r \right) \cr
&\quad+\int_s^t{\int_{\mathbb{Z}}{U_{1}\left( t,r \right)G_1\left(r, u_{\epsilon}\left( r \right) ,z \right)}}\tilde{N}_1\left( dr,dz \right), \\
v_{\epsilon}\left( t \right) &=U_{\alpha ,\epsilon,2}\left( t,s \right) y+\psi _{\alpha ,\epsilon,2}\left( v_{\epsilon};s \right)\left( t \right)  +\frac{1}{\epsilon}\int_s^t{U_{\alpha ,\epsilon,2}\left( t,r \right) B_2\left( r,u_{\epsilon}\left( r \right) ,v_{\epsilon}\left( r \right) \right)}dr \cr
&\quad +\frac{1}{\sqrt{\epsilon}}\int_s^t{U_{\alpha ,\epsilon,2}\left( t,r \right) F_2\left( r,v_{\epsilon}\left( r \right)  \right)}dw^{Q_2}\left( r \right) \cr
&\quad +\int_s^t{\int_{\mathbb{Z}}{U_{\alpha ,\epsilon,2}\left( t,r \right) G_2\left( r,v_{\epsilon}\left( r \right),z \right)}}\tilde{N}_{2}^{\epsilon}\left( dr,dz \right).
\end{cases} 
\end{split}
\end{eqnarray}

Under the  assumptions  (A1)-(A6), by proceeding as \cite[Lemma 5.1]{Xu2018Averaging} and \cite[Lemma 3.1]{cerrai2011averaging}, we can get that for any $p\ge 1$ and $T>0$, there exists a positive constant $c_{p,T}$, such that for any $x,y\in \mathbb{E}$ and $\epsilon \in \left( 0,1 \right]$, we have   
\begin{eqnarray}\label{en33}
\mathbf{E}\underset{t\in \left[ s,T \right]}{\sup}\left\| u_{\epsilon}\left( t \right) \right\| _{ \mathbb{E}}^{p}\le c_{p,T}\left( 1+\left\| x \right\| _{ \mathbb{E}}^{p}+\left\| y \right\| _{ \mathbb{E}}^{p} \right),
\end{eqnarray}
and
\begin{eqnarray}\label{en34}
\int_s^T{\mathbf{E}\left\| v_{\epsilon}\left( t \right) \right\| _{ \mathbb{E}}^{p}}dt\le c_{p,T}\left( 1+\left\| x \right\| _{ \mathbb{E}}^{p}+\left\| y \right\| _{ \mathbb{E}}^{p} \right).
\end{eqnarray}

Then, due to the equation (\ref{en5}) and the estimates (\ref{en33}) and (\ref{en34}), using the proof of \cite[Proposition 3.2]{cerrai2011averaging} to the present situation,  we can prove that there exists $ \bar{\theta}>0 $, such that for any $T>s, x\in C^\theta( \bar{\mathcal{O}}) $ with $\theta \in [ 0,\bar{\theta} )$ and $y\in E$, we have
\begin{eqnarray}\label{en313}
\underset{\epsilon \in \left( 0,1 \right]}{\sup}\mathbf{E} \left\| u_{\epsilon}\left( t \right) \right\|_{L^\infty  ( s,T; C^\theta( \bar{\mathcal{O}})  ) }\le c_{ \theta, T}\left( 1+\left\| x \right\| _{C^\theta( \bar{\mathcal{O}})}^{m_1}+\left\| y \right\| _{\mathbb{E}}^{m_1} \right),
\end{eqnarray}
where $c_{ \theta, T}>0$ is a positive constant.

Finally, by proceeding as the proof of \cite[Lemma 5.3]{Xu2018Averaging}, we can show that for any $\theta \in [0,\bar{\theta}) $, there also exists $\beta \left( \theta \right) >0$, such that, for any $T>0,p\ge 1,x\in C^{\theta}(\bar{\mathcal{O}}), y\in \mathbb{E}$ and $ s,t\in [0,T] $, we have
	\begin{eqnarray}\label{en317}
	\underset{\epsilon \in \left( 0,1 \right]}{\sup}\mathbf{E}\left\| u_{\epsilon}\left( t \right) -u_{\epsilon}\left( s \right) \right\| _{\mathbb{E}}^{p}\leq c_{p, \theta, T}\big( \left|  t-s \right| ^{\beta \left( \theta \right) p}+\left|  t-s \right|\big) \left( 1+\left\| x \right\| _{C^{\theta}(\bar{\mathcal{O}})}^{m_1p}+\left\| y \right\|_\mathbb{E}^{m_1p} \right). 
	\end{eqnarray}

According to the equation (\ref{en313}) and (\ref{en317}), using the Arzel\`{a}-Ascoli theorem, we know that the family $\left\{ \mathscr{L} \left( u_{\epsilon} \right) \right\} _{\epsilon \in \left( 0,1 \right]}$ is tight.

\section{The averaged equation}\label{sec-5}
In this section, we research the fast equation with frozen slow component $x\in E$, we main prove that there also exists an evolution family of measures for this fast equation and define the averaged equation through it.

First, for any $s\in \mathbb{R}$, any frozen slow component $x\in \mathbb{E}$ and initial condition $y\in \mathbb{E}$, we introduce the following problem
\begin{eqnarray}\label{en41}
dv\left( t \right) &=&\left[ \left( A_2\left( t \right) -\alpha \right) v\left( t \right) +B_2\left( t,x,v\left( t \right) \right) \right] dt+F_2\left( t, v\left( t \right) \right) d\bar{\omega}^{Q_2}\left( t \right) \cr
&&+{\int_{\mathbb{Z}}{G_2}\left( t, v\left( t \right) ,z \right)}{\tilde{N}_{{2}^{'}}}\left( dt,dz \right), \qquad\qquad\qquad\qquad v\left( s \right)=y,
\end{eqnarray}
where 
$$
\bar{w}^{Q_2}\left( t \right) =\left\{ \begin{array}{l}
w_{1}^{Q_2}\left( t \right),  \\
w_{2}^{Q_2}\left( -t \right),  \\
\end{array} \right. \begin{array}{c}
if\ t\ge 0,\\
if\ t<0,\\
\end{array}
$$
$$
{\tilde{N}_{{2}^{'}}}\left( t,z \right) =\left\{ \begin{array}{l}
{\tilde{N}_{{1}^{'}}}\left( t,z \right), \\
{\tilde{N}_{{3}^{'}}}\left( -t,z \right),  \\
\end{array} \right. \begin{array}{c}
if\ t\ge 0,\\
if\ t<0,\\
\end{array}
$$
where ${\tilde{N}_{{1}^{'}}}\left( dt,dz \right)$ and ${\tilde{N}_{{3}^{'}}}\left( dt,dz \right)$ has the same  L\'{e}vy measure. The process $w_{1}^{Q_2}\left( t \right)$, $w_{2}^{Q_2}\left( t \right)$, ${\tilde{N}_{{1}^{'}}}\left( dt,dz \right)$ and  ${\tilde{N}_{{3}^{'}}}\left( dt,dz \right)$ are independent and the definition of which is given in Section \ref{sec-2}.

According to the prove in Section \ref{sec-3}, we can get that for any $ x,y\in\mathbb{E},p\geq 1$ and $ s<T $, there exists a unique mild solution $ v^x\left( \cdot;s,y\right)\in  L^p\left( \Omega ;\mathbb{D}\left( \left[  s,T \right] ;\mathbb{E} \right) \right)  $. And using the same argument as  our previous work \cite{Xu2018Averaging}, we can get that there also exists $ \delta >0 $,  such that for any $ x,y\in\mathbb{E}$ and $ p\ge 1 $, we have 
\begin{eqnarray}\label{en413}
\mathbf{E}\left\| v^x\left( t;s,y \right) \right\| _{\mathbb{E}}^{p}\le c_p\left( 1+\left\| x \right\| _{\mathbb{E}}^{p}+e^{-\delta p\left( t-s \right)}\left\| y \right\| _{\mathbb{E}}^{p} \right), \quad s<t.
\end{eqnarray}

Next, same as our previous work \cite{Xu2018Averaging}, if $ t\in\mathbb{R} $, we also giving the following problem 
\begin{eqnarray}\label{en42}
dv\left( t \right) &=&\left[ \left( A_2\left( t \right) -\alpha \right) v\left( t \right) +B_2\left( t,x,v\left( t \right) \right) \right] dt+F_2\left( t, v\left( t \right) \right) d\bar{\omega}^{Q_2}\left( t \right) \cr
&&+\int_{\mathbb{Z}}{G_2}\left( t, v\left( t \right) ,z \right) \tilde{N}_{{2}^{'}}\left( dt,dz \right),
\end{eqnarray}
 for every $ s<t $.
 
 By proceeding as \cite{cerrai2017averaging} and using the conclusion we have proved in \cite{Xu2018Averaging}, it is easy to prove that 
 for any $ t\in \mathbb{R} $  and  $p\ge 1  $,  there exists $ \eta ^x\left( t \right) \in L^p\left( \varOmega ;\mathbb{E} \right) $ such that for all $ x,y\in\mathbb{E}$, we have
 \begin{eqnarray}\label{en414}
 \underset{s\rightarrow -\infty}{\lim}\mathbf{E}\left\| v^x\left( t;s,y \right) -\eta ^x\left( t \right) \right\| _{\mathbb{E}}^{p}=0.
 \end{eqnarray}
and we can get that  $ \eta ^x $ is a mild solution in $ \mathbb{R} $ of equation (\ref{en42}). Moreover, for any $ R>0 $, there also exists $  c_R>0 $ such that for any $ x_1,x_2\in\mathbb{E}$
\begin{eqnarray}\label{en0431}
x_1,x_2\in B_\mathbb{E}(R)\Rightarrow\underset{t\in \mathbb{R}}{\sup}\mathbf{E}\left\| \eta ^{x_1}(t)-\eta ^{x_2}(t) \right\| _{\mathbb{E}}^{2}\le c_R\left\| x_1-x_2 \right\| _{\mathbb{E}}^{2}.
\end{eqnarray}

Then, for any  $ t\in \mathbb{R} $ and  $ x\in\mathbb{E}$, we denote that the law of the random variable  $ \eta ^x\left( t \right)  $ is $ \mu _{t}^x  $. As the prove of our previous work \cite{Xu2018Averaging}, we also can get that $\left\{ \mu _{t}^x  \right\}_{t \in \mathbb{R} }$ defines an evolution family of measures on $\mathbb{E}$ for equation (\ref{en41}). 

Now, we give the following assumption.
\begin{enumerate}[({A}7)]
	\item 
	\begin{enumerate}
		\item The functions $ \gamma_2:\mathbb{R}\rightarrow \left( 0,\infty \right)  $  and $ l_2:\mathbb{R}\times \mathcal{O}\rightarrow \mathbb{R}^d $  are periodic, with the same period.
		\item The families of functions
		\begin{gather}
		\mathbf{B}_{1,R}:=\left\{ b_1\left( \cdot ,\xi ,\sigma \right) :\ \xi \in \mathcal{O},\ \sigma \in \mathcal{B}_{\mathbb{R}^2}\left( R \right) \right\},\cr
		\mathbf{B}_{2,R}:=\left\{ b_2\left( \cdot ,\xi ,\sigma \right) :\ \xi \in \mathcal{O},\ \sigma \in \mathcal{B}_{\mathbb{R}^2}\left( R \right) \right\},\cr
		\mathbf{F}_R:=\left\{ f_2\left( \cdot ,\xi ,\sigma \right) :\ \xi \in \mathcal{O},\ \sigma \in \mathcal{B}_{\mathbb{R} }\left( R \right) \right\},\cr
		\mathbf{G}_R:=\left\{ g_2\left( \cdot ,\xi ,\sigma ,z \right) :\ \xi \in \mathcal{O},\ \sigma \in \mathcal{B}_{\mathbb{R} }\left( R \right) ,\ z\in \mathbb{Z} \right\},\nonumber
		\end{gather}
		are uniformly almost periodic for any  $ R>0 $. 
	\end{enumerate}
\end{enumerate}
\begin{rem}\label{rem4.5}
	{\rm Similar with the proof of \cite[Lemma 6.2]{cerrai2017averaging}, we get that under the assumption (A7), for any  $ R>0 $, the families of functions
		\begin{gather}
		\left\{ B_1\left( \cdot ,x,y \right) :\left( x,y \right) \in \mathcal{B}_{\mathbb{E}\times \mathbb{E}}\left( R \right) \right\}, \quad
		\left\{ B_2\left( \cdot ,x,y \right) :\left( x,y \right) \in \mathcal{B}_{\mathbb{E}\times \mathbb{E}}\left( R \right) \right\}, \cr
		\left\{ F_2\left( \cdot , y \right) : y \in \mathcal{B}_{\mathbb{E} }\left( R \right) \right\},\quad
		\left\{ G_2\left( \cdot , y,z \right) :\left( y,z \right) \in \mathcal{B}_{\mathbb{E} }\left( R \right) \times \mathbb{Z} \right\}, \nonumber
		\end{gather}
		are uniformly almost periodic.}
\end{rem}

As in \cite{Xu2018Averaging} and \cite{cerrai2017averaging}, we can prove that under the assumptions (A1)-(A7), the mapping $t\in \mathbb{R}\mapsto \mu ^x_{t} \in \mathcal{P}\left(\mathbb{E}\right)$ is almost periodic. Then, due to (\ref{en0431}), we also can get that for every compact set  $ \mathbb{K} \subset \mathbb{E}$, the family of functions
\begin{eqnarray}\label{en61}
\left\{ t\in \mathbb{R}\mapsto \int_{ \mathbb{E}}{B_1\left(t, x,y \right)}\mu ^x_{t} \left( dy \right):\quad x\in \mathbb{K} \right\} 
\end{eqnarray}
is uniformly almost periodic. So, we define 
\begin{eqnarray}\label{en17}
\bar{B}_1\left( x \right) :=\underset{T\rightarrow \infty}{\lim}\frac{1}{T}\int_0^T{\int_{ \mathbb{E}}{B_1\left( t,x,y \right)}\mu _{t} \left( dy \right) dt, \quad x\in  \mathbb{E}},
\end{eqnarray}
we can get that the mapping $ \bar{B}_1:\mathbb{E}\rightarrow\mathbb{E}$ is locally Lipschitz-continuous. Similar with the prove of \cite[Lemma 4.2]{Xu2018Averaging} and \cite[Lemma 7.2]{cerrai2017averaging}, we can conclude that the following crucial results are also established in this paper.
\begin{lem}\label{lem6.2}Under the assumptions {\rm (A1)-(A7)}, for any  $ T>0,s\in \mathbb{R} $ and  $ x,y\in \mathbb{E}$, there exist some constants $ \kappa_1,\kappa_2\geq 0 $, we have 
	\begin{eqnarray}\label{en63}
	\mathbf{E}\Big| \frac{1}{T}\int_t^{t+T}{  B_1\left(t, x,v^x\left( t;s,y \right) \right)   }dr-  \bar{B}_1\left( x \right) \Big|  
	&\le& \frac{c}{T }\left( 1+\left\| x \right\| _{ \mathbb{E}}^{\kappa_1} +\left\| y \right\| _{ \mathbb{E}}^{\kappa_2}  \right) +\alpha \left( T,x \right) 	
	\end{eqnarray}
	for some mapping  $ \alpha : \left[ 0,\infty \right) \times  \mathbb{E}\rightarrow \left[ 0,\infty \right)  $ such that
	\begin{eqnarray}\label{en31}
	\underset{T>0}{\sup}\alpha \left( T,x \right) \le c\left( 1+\left\| x \right\| _{ \mathbb{E}}^{m_1}  \right), \quad x\in  \mathbb{E},
	\end{eqnarray}
	and for any compact set $ \mathbb{K}\subset\mathbb{E}$, have 
	\begin{eqnarray}\label{en32}
	\underset{T\rightarrow \infty}{\lim}\underset{x\in  \mathbb{E}}{\sup}\ \alpha \left( T,x \right) =0.
	\end{eqnarray}	
\end{lem}

We introduce the following averaged equation
\begin{eqnarray}\label{en67}
du\left( t \right) &=&\left[ A_1(t)u\left( t \right) +\bar{B}_1\left( u\left( t \right) \right) \right] dt+F_1\left(t, u\left( t \right) \right) dw^{Q_1}\left( t \right) \cr
&&+\int_{\mathbb{Z}}{G_1\left( t,u\left( t \right) ,z \right)}\tilde{N}_1\left( dt,dz \right), \quad\quad  u\left( 0 \right) =x\in  \mathbb{E}.
\end{eqnarray}
Due to Theorem \ref{th3.4}, we can prove that for any  $ T>0, p\ge 1 $ and  $ x\in \mathbb{E}$, equation  (\ref{en67}) admits a unique mild solution $\bar{u}$.

\section{Averaging principles}\label{sec-6}
In this section, we will show that the validity of the averaging principle. That is, the slow motion $ u_{\epsilon}  $ will converges to the averaged motion  $ \bar{u} $, as $ \epsilon \rightarrow 0 $.
\begin{thm}\label{thm7.1}Under the assumptions (A1)-(A7), fix $ x\in C ^{\theta}(\bar{\mathcal{O}}) $ with $ \theta \in [ 0,\bar{\theta} ) $, and  $ y\in\mathbb{E}$, for any  $ T>0 $ and  $ \eta >0 $, we have 
	\begin{eqnarray}\label{en71}
	\underset{\epsilon \rightarrow 0}{\lim}\mathbb{P}\Big( \underset{t\in \left[ 0,T \right]}{\sup}\left\| u_{\epsilon}\left( t \right) -\bar{u}\left( t \right) \right\| _{\mathbb{E}}>\eta \Big) =0,
	\end{eqnarray}
	where  $ \bar{u} $ is the solution of the averaged equation (\ref{en67}).
\end{thm}
\para{Proof:} For any  $ h\in \mathcal{D}\left( A_1 \right) $, we have
\begin{eqnarray}
\int_\mathcal{O}{u_{\epsilon}\left( t,\xi \right)}h\left( \xi \right) d\xi &=&\int_\mathcal{O}{x\left( \xi \right)}h\left( \xi \right) d\xi +\int_0^t{\int_\mathcal{O}{u_{\epsilon}\left( r,\xi \right)}A_1\left( r \right) h\left( \xi \right)}d\xi dr\cr
&&+\int_0^t{\int_\mathcal{O}{\bar{B}_1\left( u_{\epsilon}\left( r \right) \right) \left( \xi \right) h\left( \xi \right)}}d\xi dr+\int_0^t{\int_\mathcal{O}{\left[ F_1\left( r,u_{\epsilon}\left( r \right) \right) h \right] \left( \xi \right)}}dw^{Q_1}\left( r,\xi \right) \cr
&&+\int_0^t{\int_{\mathbb{Z}}{\int_\mathcal{O}{\left[ G_1\left( r,u_{\epsilon}\left( r \right) ,z \right) h \right]}}\left( \xi \right)}d\xi\tilde{N}_1\left(  dr,dz \right)+R_{\epsilon}\left( t \right) ,\nonumber
\end{eqnarray} 
where
$$
R_{\epsilon}\left( t \right) :=\int_0^t{\int_\mathcal{O}{\left( B_1\left( r,u_{\epsilon}\left( r \right) ,v_{\epsilon}\left( r  \right) \right) \left( \xi \right) -\bar{B}_1\left( u_{\epsilon}\left( r  \right) \right) \left( \xi \right) \right) h\left( \xi \right)}}d\xi dr.
$$

As the proof in \cite{Xu2018Averaging}, because we have get that the family $ \left\{ \mathscr{L}\left( u_{\epsilon} \right) \right\} _{\epsilon \in \left( 0,1 \right]} $  is tight in Section \ref{sec-4}. If we want to prove  \thmref{thm7.1}, it is sufficient to prove  that for any $ T> 0$, we have $\underset{\epsilon \rightarrow 0}{\lim}\mathbf{E}\underset{t\in \left[ 0,T \right]}{\sup}\left\| R_{\epsilon}\left( t \right) \right\|_\mathbb{E}=0 $.

First, for any $ n\in \mathbb{N} $, we define 
\begin{eqnarray}
b_{i,n}\left( t,\xi ,\sigma _1,\sigma _2 \right) :=\left\{ \begin{array}{c}
b_i\left( t,\xi ,\sigma _1,\sigma _2 \right) ,\\
b_i\left( t,\xi ,\sigma _1n/\left| \sigma _1 \right|,\sigma _2 \right) ,\\
\end{array} \right. \quad \begin{array}{c}
if\ \left| \sigma _1 \right|\leq n,\\
if\ \left| \sigma _1 \right|>n.\\
\end{array}
\end{eqnarray}
For each $ b_{i,n} $, denote the corresponding composition operator is $ B_{i,n} $, and we have
\begin{eqnarray}\label{en061}
x\in B_\mathbb{E}\left( n \right) \Rightarrow B_{i,n}\left( t,x,y \right) =B_i\left( t,x,y \right) ,\ \quad t\in \mathbb{R},\  y\in \mathbb{E}.
\end{eqnarray}
It is easy to get that the mapping $ b_{1,n} $ and $ b_{2,n} $ satisfy all conditions in (A3) and (A4), respectively. And for any fixed $ \left( t,\xi\right) \in  \mathbb{R} \times \bar{\mathcal{O}} $ and $ \sigma_2\in \mathbb{R} $, the mapping $ b_{i,n}\left( t,\xi ,\cdot,\sigma _2 \right) $ are Lipschitz-continuous. 

In addition, for any $ n\in\mathbb{N} $, we define
\begin{eqnarray}
f_{1,n}\left( t,\xi ,\sigma \right) :=\left\{ \begin{array}{c}
f_1\left( t,\xi ,\sigma \right) ,\\
f_1\left( t,\xi ,\sigma n/\left| \sigma \right| \right) ,\\
\end{array} \right. \quad \begin{array}{c}
if\ \left| \sigma \right|\leq n,\\
if\ \left| \sigma \right|>n.\\
\end{array},\nonumber
\end{eqnarray}
and
\begin{eqnarray}
g_{1,n}\left( t,\xi ,\sigma \right) :=\left\{ \begin{array}{c}
g_1\left( t,\xi ,\sigma ,z \right) ,\\
g_1\left( t,\xi ,\sigma n/\left| \sigma \right|,z \right) ,\\
\end{array} \right. \quad \begin{array}{c}
if\ \left| \sigma \right|\leq n,\\
if\ \left| \sigma \right|>n.\\
\end{array},\nonumber
\end{eqnarray}
where $ (t,\xi)\in\mathbb{R}\times\bar{\mathcal{O}} $ and $ z\in \mathbb{Z} $. The corresponding composition operator of $ f_{1,n} $ and $ g_{1,n} $ are denoted by $ F_{1,n} $ and $ G_{1,n} $, respectively.

Now, for any $ n\in\mathbb{N} $, we introduce the following system
\begin{eqnarray}\label{en0611}
\begin{split}
\begin{cases}
du \left( t \right) &=\left[ A_1\left( t\right) u \left( t \right) +B_{1,n}\left(t, u \left( t \right) ,v \left( t \right) \right) \right] dt+F_{1,n}\left(t, u \left( t \right) \right) d\omega ^{Q_1}\left( t \right) \\
&\quad +\int_{\mathbb{Z}}{G_{1,n}\left(t, u \left( t \right) ,z \right)}\tilde{N}_1\left( dt,dz \right), \\
dv \left( t \right) &=\frac{1}{\epsilon}\left[ \left( A_2\left( t \right) -\alpha \right) v \left( t \right) +B_{2,n}\left( t,u \left( t \right) ,v \left( t \right) \right) \right] dt \\
&\quad +\frac{1}{\sqrt{\epsilon}}F_2\left( t,v \left( t \right) \right) d\omega ^{Q_2}\left( t \right)  +\int_{\mathbb{Z}}{G_2}\left( t,v \left( t \right)  ,z\right) \tilde{N}_{2}^{\epsilon}\left( dt,dz \right), \\
u \left( s \right) &=x, \quad v \left( s \right) =y,
\end{cases}
\end{split}
\end{eqnarray}
we denote the solution of (\ref{en0611}) is $ \left( u_{\epsilon,n},v_{\epsilon,n}\right)  $.

Then, for any $ n\in \mathbb{N} $  and any frozen slow component $x\in \mathbb{E}$, we introduce the following problem
\begin{eqnarray}\label{en062}
dv\left( t \right) &=&\left[ \left( A_2\left( t \right) -\alpha \right) v\left( t \right) +B_{2,n}\left( t,x,v\left( t \right) \right) \right] dt+F_2\left( t,v\left( t \right) \right) d\omega ^{Q_2}\left( t \right)\cr
&& +\int_{\mathbb{Z}}{G_2}\left( t,v\left( t \right) ,z \right) \tilde{N}_2\left( dt,dz \right) ,\qquad\qquad\qquad\qquad v(s)=y,                         
\end{eqnarray}
and denote its solution is $ v_n^x(t;s,y) $. Thanks to (\ref{en061}), for any $ t\geq0 $ and $ x\in\mathbb{E}$, we have
\begin{eqnarray}
v_{n}^{x}\left( t;s,y \right) =\left\{ \begin{array}{c}
v^x\left( t;s,y \right) ,\\
v^{x_n}\left( t;s,y \right) ,\\
\end{array} \right. \quad \begin{array}{c}
if\ \left| x\left( \xi \right) \right|\leq n,\\
if\ \left| x\left( \xi \right) \right|>n,\\
\end{array}\nonumber
\end{eqnarray}
where $ x_n=nsignx(\xi) $.

Due to the coefficients of equation (\ref{en062}) satisfy the same conditions as the equation (\ref{en42}), for each $ x\in\mathbb{E}$, there exists an evolution of measures family $ \left\lbrace \mu_t^{x,n}\right\rbrace_{t\in\mathbb{R}} $ for equation (\ref{en062})
\begin{eqnarray}
\mu _{t}^{x,n}=\left\{ \begin{array}{c}
\mu _{t}^{x},\\
\mu _{t}^{x_n},\\
\end{array} \right. \quad \begin{array}{c}
if\ \left| x\left( \xi \right) \right|\leq n,\\
if\ \left| x\left( \xi \right) \right|>n.\\
\end{array}\nonumber
\end{eqnarray}
As proof of (\ref{en413}), for any $ T>0, x\in\mathbb{E}$ and $ p\geq 1 $, we can get that there also exists $ \delta>0 $, such that
\begin{eqnarray}\label{en064}
\mathbf{E}\left\| v_n^x\left( t;s,y \right) \right\| _{\mathbb{E}}^{p}\le c_{p,n}\left( 1+ e^{-\delta p\left( t-s \right)}\left\| y \right\| _{\mathbb{E}}^{p} \right), \quad s<t.
\end{eqnarray}	
 Similarly, we can define 
\begin{eqnarray}\label{en063}
\bar{B}_{1,n}\left( x \right) :=\underset{T\rightarrow \infty}{\lim}\frac{1}{T}\int_0^T{\int_{ E}{B_{1,n}\left( t,x,y \right)}\mu _{t}^{x,n}\left( dy \right) dt, \quad x\in  E},
\end{eqnarray}
 and, we have 
$$
\left\| x \right\|_\mathbb{E}\le n\Rightarrow \bar{B}_{1,n}\left( x \right) =\bar{B}_1\left( x \right).
$$
Moreover, it is easy to prove that the mapping $ \bar{B}_{1,n}:\mathbb{E}\rightarrow \mathbb{E}$ is Lipschitz-continuous, and the results similar with (\ref{en63})-(\ref{en32})  also can be established.

Next, we prove that the validity of the averaging principle by using the classical Khasminskii method as in \cite{Xu2018Averaging}. For any  $ \epsilon >0 $, we divide the interval $ \left[ 0,T \right]  $ in subintervals of the size  $ \delta _{\epsilon} $, where $ \delta _{\epsilon}>0 $ is a deterministic constant. Then, we define the following auxiliary fast motion $ \hat{v}_{\epsilon} $ in each time interval  $ \left[ k\delta _{\epsilon},\left( k+1 \right) \delta _{\epsilon} \right] ,k=0,1,\cdots ,\lfloor  {T}/{\delta _{\epsilon}} \rfloor  $ 
\begin{eqnarray}\label{en72}
\begin{split}
\begin{cases}
d\hat{v}_{\epsilon ,n}\left( t \right) &=\frac{1}{\epsilon}\left[ \left( A_2\left( t \right) -\alpha \right) \hat{v}_{\epsilon ,n}\left( t \right) +B_{2,n}\left( t,u_{\epsilon ,n}\left( k\delta _{\epsilon} \right) ,\hat{v}_{\epsilon ,n}\left( t \right) \right) \right] dt\cr
&\quad+\frac{1}{\sqrt{\epsilon}}F_2\left( t,\hat{v}_{\epsilon ,n}\left( t \right) \right) d\omega ^{Q_2}\left( t \right)+\int_{\mathbb{Z}}{G_2}\left( t,\hat{v}_{\epsilon ,n}\left( t \right), z\right) \tilde{N}_{2}^{\epsilon}\left( dt,dz \right) ,\cr
\hat{v}_{\epsilon ,n}\left( k\delta _{\epsilon} \right) &=v_{\epsilon ,n}\left( k\delta _{\epsilon} \right) .
\end{cases}
\end{split}
\end{eqnarray}
Like the equation (\ref{en34}), we also can prove that for any  $ p\ge 1 $, we have 
\begin{eqnarray}\label{en73}
\int_0^T{\mathbf{E}\left\| \hat{v}_{\epsilon,n}\left( t \right) \right\| _{\mathbb{E}}^{p}}dt\le c_{p,T}\left( 1+\left\| x \right\| _{\mathbb{E}}^{p}+\left\| y \right\| _{\mathbb{E}}^{p} \right).
\end{eqnarray}
\begin{lem}\label{lem7.2}Under the assumptions (A1)-(A7), fix $ x\in C ^{\theta}(\bar{\mathcal{O}}) $ with $ \theta \in [ 0,\bar{\theta} ) $, and  $ y\in\mathbb{E}$, there exists a constant  $ \kappa >0 $, such that if
	$$
	\delta _{\epsilon}=\epsilon \ln ^{\epsilon ^{-\kappa}},
	$$
	and, for any fixed $ n\in \mathbb{N} $, we have 
	\begin{eqnarray}\label{en74}
	\underset{\epsilon \rightarrow 0}{\lim}\underset{t\in \left[ 0,T \right]}{\sup}\mathbf{E}\left\| \hat{v}_{\epsilon,n}\left( t \right) -v_{\epsilon,n}\left( t \right) \right\| _{\mathbb{E}}^{p}=0.
	\end{eqnarray}
\end{lem}
\para{Proof:} Fixed $ \epsilon >0 $ and $ n\in \mathbb{N} $. For any $ t\in \left[ k\delta _{\epsilon},\left( k+1 \right) \delta _{\epsilon} \right], k=0,1,\cdots ,\lfloor {T}/{\delta _{\epsilon}} \rfloor $, let $ \rho_{\epsilon,n}\left( t \right)  $  be the solution of the following problem
\begin{eqnarray}
d\rho _{\epsilon,n}\left( t \right) &=&\frac{1}{\epsilon}\left( A_2\left( t \right) -\alpha \right) \rho _{\epsilon,n}\left( t \right) dt+\frac{1}{\sqrt{\epsilon}}K_{\epsilon,n}\left( t \right) d\omega ^{Q_2}\left( t \right) \cr
&&+\int_{\mathbb{Z}}{H_{\epsilon,n}\left( t,z \right)}\tilde{N}_{2}^{\epsilon}\left( dt,dz \right), \quad\quad\qquad \rho_{\epsilon}\left( k\delta _{\epsilon} \right) =0,\nonumber
\end{eqnarray} 
where
$$
K_{\epsilon,n}\left( t \right) :=F_2\left( t, \hat{v}_{\epsilon,n}\left( t \right) \right) -F_2\left( t, v_{\epsilon,n}\left( t \right) \right), 
$$
$$
H_{\epsilon,n}\left( t,z \right) :=G_2\left( t, \hat{v}_{\epsilon,n}\left( t \right) ,z \right) -G_2\left( t, v_{\epsilon,n}\left( t \right) ,z \right). 
$$
We have 
$$
\rho _{\epsilon,n}\left( t \right) =\psi_{\alpha ,\epsilon,2}\left( \rho_{\epsilon,n};k\delta _{\epsilon} \right) \left( t \right) +\varGamma _{\epsilon,n}\left( t \right) +\varPsi _{\epsilon,n}\left( t \right), \quad t\in \left[ k\delta _{\epsilon},\left( k+1 \right) \delta _{\epsilon} \right], 
$$
where
$$
\varGamma _{\epsilon,n}\left( t \right) =\frac{1}{\sqrt{\epsilon}}\int_{k\delta _{\epsilon}}^t{U_{\alpha ,\epsilon,2}\left( t,r \right) K_{\epsilon,n}\left( r \right)}dw^{Q_2}\left( r \right), 
$$
$$
\varPsi _{\epsilon,n}\left( t \right) =\int_{k\delta _{\epsilon}}^t{\int_{\mathbb{Z}}{U_{\alpha ,\epsilon,2}\left( t,r \right) H_{\epsilon,n}\left( r,z \right)}}\tilde{N}_{2}^{\epsilon}\left( dr,dz \right). 
$$
Using the same arguement as \cite[Lemma 6.3]{cerrai2011averaging}, we yield 
\begin{eqnarray}\label{en202}
\mathbf{E}\left\| \varGamma _{\epsilon,n}\left( t \right) \right\| _{\mathbb{E}}^{p} \leq \frac{c_{p,n}}{\epsilon}\int_{k\delta _{\epsilon}}^{t}{\mathbf{E}\left\| \hat{v}_{\epsilon,n}\left( r \right) -v_{\epsilon,n}\left( r \right) \right\| _{\mathbb{E}}^{p} }dr.
\end{eqnarray}
For $ \varPsi _{\epsilon,n}\left( t \right) $, using  Kunita's first inequality and the H\"{o}lder inequality, we can get  
\begin{eqnarray}\label{en203}
\mathbf{E}\left\| \varPsi _{\epsilon,n}\left( t \right) \right\| _{   \mathbb{E}}^{p}
&\leq& c_{p,n}\mathbf{E}\Big(\frac{1}{\epsilon}\int_{k\delta _{\epsilon}}^t{\int_{\mathbb{Z}}{\left\| e^{-\frac{\alpha}{\epsilon}\left( t-r \right)}e^{\frac{\gamma_2\left(t,r \right)  }{\epsilon}A_2}H_{\epsilon,n}\left( r,z \right)\right\| } _{  _\mathbb{E}}^{2}}v_2\left( dz \right) dr \Big)^{\frac{p}{2}}\cr
&&+\frac{c_{p,n}}{\epsilon}\mathbf{E}\int_{k\delta _{\epsilon}}^t{\int_{\mathbb{Z}}{\left\| e^{-\frac{\alpha}{\epsilon}\left( t-r \right)}e^{\frac{\gamma_2\left(t,r \right)  }{\epsilon}A_2}H_{\epsilon,n}\left( r,z \right)\right\| } _{  _\mathbb{E}}^{p}}v_2\left( dz \right) dr \cr
&\leq& \frac{c_{p,n}}{\epsilon^{\frac{p}{2}}}\mathbf{E}\Big(\int_{k\delta _{\epsilon}}^t{\int_{\mathbb{Z}}{\left\|  H_{\epsilon,n}\left( r,z \right)\right\| } _{  _\mathbb{E}}^{2}}v_2\left( dz \right)dr\Big)^{\frac{p}{2}}\cr
&&+\frac{c_{p,n}}{\epsilon}\mathbf{E}\int_{k\delta _{\epsilon}}^t{\int_{\mathbb{Z}}{\left\|  H_{\epsilon,n}\left( r,z \right)\right\| } _{  _\mathbb{E}}^{p}}v_2\left( dz \right) dr \cr
&\leq& c_{p,n}\big( {\delta_{\epsilon}^{\frac{p-2}{2}}}/{\epsilon^{\frac{p}{2}}}+{1}/{\epsilon} \big) \int_{k\delta _{\epsilon}}^{t}{\mathbf{E}\left\| \hat{v}_{\epsilon,n}\left( r \right) -v_{\epsilon,n}\left( r \right) \right\| _{  _\mathbb{E}}^{p} }dr.
\end{eqnarray}
Thanks to $ \alpha>0 $ is large enough, we have
\begin{eqnarray}\label{en12}
\lVert \rho _{\epsilon ,n}\left( t \right) \rVert _\mathbb{E}&\leq& \lVert \psi _{\alpha ,\epsilon ,2}\left( \rho _{\epsilon ,n};k\delta _{\epsilon} \right) \left( t \right) \rVert_\mathbb{E}+\lVert \varGamma _{\epsilon ,n}\left( t \right) \rVert_\mathbb{E}+\lVert \varPsi _{\epsilon ,n}\left( t \right) \rVert_\mathbb{E}\cr
&\leq& c_{p,n}\big( {\delta_{\epsilon}^{\frac{p-2}{2}}}/{\epsilon^{\frac{p}{2}}}+{1}/{\epsilon} \big) \int_{k\delta _{\epsilon}}^{t}{\mathbf{E}\left\| \hat{v}_{\epsilon,n}\left( r \right) -v_{\epsilon,n}\left( r \right) \right\| _{   \mathbb{E}}^{p} }dr.
\end{eqnarray}
If we denote $ \varLambda _{\epsilon,n}\left( t \right) :=\hat{v}_{\epsilon,n}\left( t \right) -v_{\epsilon,n}\left( t \right)  $  and $ \vartheta _{\epsilon,n}\left( t \right) := \varLambda _{\epsilon,n}\left( t \right) -  \rho _{\epsilon,n}\left( t \right)$, we have
\begin{eqnarray}
d\vartheta _{\epsilon,n}\left( t \right) &=&\frac{1}{\epsilon}\left[ \left( A_2\left( t \right) -\alpha \right) \vartheta _{\epsilon,n}\left( t \right) +B_{2,n}\left( t,u_{\epsilon,n}\left( k\delta _{\epsilon} \right) ,\hat{v}_{\epsilon,n}\left( t \right) \right) -B_{2,n}\left( t,u_{\epsilon,n}\left( t \right) ,v_{\epsilon,n}\left( t \right) \right) \right] dt\cr
&=&\frac{1}{\epsilon}\left[ \left( A_2\left( t \right) -\alpha \right) \vartheta _{\epsilon ,n}\left( t \right) +B_{2,n}\left( t,u_{\epsilon ,n}\left( k\delta _{\epsilon} \right) ,\hat{v}_{\epsilon ,n}\left( t \right) \right) -B_{2,n}\left( t,u_{\epsilon ,n}\left( t \right) ,\hat{v}_{\epsilon ,n}\left( t \right) \right) \right. \cr
&&\left. -\tau \left( t,u_{\epsilon ,n}\left( t \right) ,\hat{v}_{\epsilon ,n}\left( t \right) ,v_{\epsilon ,n}\left( t \right) \right) \left( \vartheta _{\epsilon ,n}\left( t \right) +\rho _{\epsilon ,n}\left( t \right) \right) \right] dt
\end{eqnarray}
By proceeding as \cite[Lemma 6.3]{cerrai2011averaging}, we can get
\begin{eqnarray}\label{en13}
\lVert \vartheta _{\epsilon ,n}\left( t \right) \rVert_\mathbb{E}&\le& \frac{c_n}{\epsilon}\int_{k\delta _{\epsilon}}^{t}{ e^{-\frac{\alpha}{\epsilon}\left( t-r \right)}\lVert u_{\epsilon ,n}\left( k\delta _{\epsilon} \right) -u_{\epsilon ,n}\left( r \right) \rVert_\mathbb{E}
}dr \cr
&&+\frac{1}{\epsilon}\underset{r\in \left[ k\delta _{\epsilon},t \right]}{\text{sup}}\lVert \rho _{\epsilon ,n}\left( r \right) \rVert_\mathbb{E}\int_{k\delta _{\epsilon}}^t{\exp \Big( -\frac{1}{\epsilon}\int_r^t{\tau _{\epsilon ,n}\left( \sigma \right)}d\sigma \Big) \tau _{\epsilon ,n}\left( r \right)}dr\cr
&\leq& c_{p,n} \big( 1+\left\| x \right\| _{C^{\theta}(\bar{\mathcal{O}})}^{m_1p}+\left\| y \right\|_\mathbb{E}^{p} \big)\big( \delta_{\epsilon} ^{\beta \left( \theta \right) p}+\delta_\epsilon\big)\cr
&&+c_{p,n}\big( {\delta_{\epsilon}^{\frac{p-2}{2}}}/{\epsilon^{\frac{p}{2}}}+{1}/{\epsilon} \big) \int_{k\delta _{\epsilon}}^{t}{\mathbf{E}\left\| \hat{v}_{\epsilon,n}\left( r \right) -v_{\epsilon,n}\left( r \right) \right\| _{   \mathbb{E}}^{p} }dr.
\end{eqnarray}
where 
$$
\tau _{\epsilon ,n}\left( r \right) :=\tau \left( \xi _{\epsilon ,n}\left( t \right) ,u_{\epsilon ,n}\left( t,\xi _{\epsilon ,n}\left( t \right) \right) ,\hat{v}_{\epsilon ,n}\left( t,\xi _{\epsilon ,n}\left( t \right) \right) ,v_{\epsilon ,n}\left( t,\xi _{\epsilon ,n}\left( t \right) \right) \right) ,
$$
and $ \xi _{\epsilon ,n}\left( t \right) \in \bar{\mathcal{O}} $, satisfy
$$
|\vartheta _{\epsilon ,n}\left( t,\xi _{\epsilon ,n}\left( t \right) \right) |=\lVert \vartheta _{\epsilon ,n}\left( t \right) \rVert _\mathbb{E}.
$$
Due to (\ref{en12}) and (\ref{en13}), for any $ p\geq 1 $, we have
\begin{eqnarray}
\mathbf{E}\left\| \hat{v}_{\epsilon}\left( t \right) -v_{\epsilon}\left( t \right) \right\| _{   \mathbb{E}}^{p} 
&\leq& c_{p,n} \big( 1+\left\| x \right\| _{C^{\theta}(\bar{\mathcal{O}})}^{m_1p}+\left\| y \right\|_\mathbb{E}^{m_1p} \big)\big( \delta_{\epsilon} ^{\beta \left( \theta \right) p}+\delta_\epsilon\big)\cr
&&+c_{p,n}\big( {\delta_{\epsilon}^{\frac{p-2}{2}}}/{\epsilon^{\frac{p}{2}}}+{1}/{\epsilon} \big) \int_{k\delta _{\epsilon}}^{t}{\mathbf{E}\left\| \hat{v}_{\epsilon,n}\left( r \right) -v_{\epsilon,n}\left( r \right) \right\| _{   \mathbb{E}}^{p} }dr.\nonumber
\end{eqnarray}
From the Gronwall lemma, this means
$$
\mathbf{E}\left\| \hat{v}_{\epsilon,n}\left( t \right) -v_{\epsilon,n}\left( t \right) \right\| _{\mathbb{E}}^{p}\leq c_{p,n} \big( 1+\left\| x \right\| _{C^{\theta}(\bar{\mathcal{O}})}^{m_1p}+\left\| y \right\|_\mathbb{E}^{p} \big)\big( \delta_{\epsilon} ^{\beta \left( \theta \right) p}+\delta_\epsilon\big)e^{c_{p,n} ( {\delta_{\epsilon}^{\frac{p-2}{2}}}/{\epsilon^{\frac{p}{2}}}+{1}/{\epsilon}  ) \delta _{\epsilon}}.
$$
For  $ t\in \left[ 0,T \right]  $, selecting  $ \delta _{\epsilon}=\epsilon \ln ^{\epsilon ^{-\kappa}} $, then if we take  $ \kappa <\frac{\beta \left( \theta \right)p}{\beta \left( \theta \right)p+2 c_{p,n}}\land \frac{1}{1+2c_{p,n}} $, we have (\ref{en74}).\qed

Finally, under the same assumptions as in \thmref{thm7.1}, by proceeding as \cite[Lemma 8.2]{cerrai2017averaging} and \cite[Lemma 6.4]{Xu2018Averaging}, for any  $ T>0 $, we can get 
\begin{eqnarray}\label{en75}
\underset{\epsilon \rightarrow 0}{\lim}\mathbf{E}\underset{t\in \left[ 0,T \right]}{\sup}\left\| R_{\epsilon}\left( t \right) \right\|_\mathbb{E}=0.\nonumber
\end{eqnarray} 

Through the above proof, \thmref{thm7.1} is established.\qed

\section{Conclusions}
In this paper, we study the averaging principle for a class of non-autonomous slow-fast system with polynomial growth. First, using the Sobolev embedding theorem, fxed point theorem and stopping technique, the existence and uniqueness of the mild solution is proved. Next, by means of the comparison theorem and the properties of transition operator, the existence of time-dependent evolution family of measures associated with the fast equation is studied, and the averaged coefcient is obtained. Finally, through the truncation technique, the averaging principle for a class of non-autonomous slow-fast systems with polynomial growth is presented.

\section*{Acknowledgments}
The research was supported in part by the NSF of China (11572247, 11802216) and the Seed Foundation of Innovation and Creation for Graduate Students in Northwestern Polytechnical University (ZZ2018027). B. Pei was an International Research Fellow of Japan Society for the Promotion of Science (Postdoctoral Fellowships for Research in Japan (Standard)). 
Y. Xu would like to thank the Alexander von Humboldt Foundation for the support.

\section*{References}
\bibliography{references}
\end{document}